\documentclass[aap,MSNbibl,nameyear,seceqn,dvips]{arximspdf}
\usepackage{graphicx}

\doi{10.1214/11-AAP836} %kopijuoti is PTS
\volume{22}
\issue{5}
\pubyear{2012}
\firstpage{2108}
\lastpage{2132}

\makeatletter
\newcommand{\eqref}[1]{(\ref{#1})}
\renewcommand{\citep}[1]{(\citeauthor{#1} \citeyear{#1})}

\newtheorem{theorem}{Theorem}
\newtheorem{proposition}{Proposition}[section]
\newtheorem{lemma}[proposition]{Lemma}
\newtheorem{corollary}[proposition]{Corollary}

\newproclaim{definition}[proposition]{Definition}
\newproclaim{remark}[proposition]{Remark}
\newproclaim{example}[proposition]{Example}
\makeatother

\begin{document}
\begin{frontmatter}

\title{Muller's ratchet with compensatory mutations}
\runtitle{Muller's ratchet with compensatory mutations}

\begin{aug}
\author[a]{\fnms{P.} \snm{Pfaffelhuber}\corref{}\ead[label=e1]{p.p@stochastik.uni-freiburg.de}\thanksref{t1}},
\author[a]{\fnms{P. R.} \snm{Staab}}
\and
\author[b]{\fnms{A.} \snm{Wakolbinger}}
\thankstext{t1}{Supported by the BMBF through FRISYS (Kennzeichen 0313921) and by
the DFG through Grant Pf672/1-1.}
\runauthor{P. Pfaffelhuber, P. R. Staab and A. Wakolbinger}
\affiliation{Albert-Ludwigs-Universit\"at Freiburg,
Albert-Ludwigs-Universit\"at Freiburg and~Goethe-Universit\"at~Frankfurt}
\address[a]{P. Pfaffelhuber\\
P. R. Staab\\
Fakult\"at f\"ur Mathematik und Physik\\
Albert-Ludwigs-Universit\"at Freiburg\\
Eckerstra\ss e 1\\
79104 Freiburg\\
Germany\\
\printead{e1}}

\address[b]{A. Wakolbinger\\
Fachbereich Informatik und Mathamatik\\
Goethe-Universit\"at Frankfurt\\
Robert-Mayer-Str. 10\\
60054 Frankfurt\\
Germany}
\end{aug}

% HISTORY:
\received{\smonth{8} \syear{2011}}
\revised{\smonth{11} \syear{2011}}

% ABSTRACT
%
\begin{abstract}
We consider an infinite-dimensional system of stochastic
differential equations describing the evolution of type frequencies
in a large population. The type of an individual is the number of
deleterious mutations it carries, where fitness of individuals
carrying $k$ mutations is decreased by $\alpha k$ for some $\alpha
>0$. Along the individual lines of descent, new mutations
accumulate at rate $\lambda$ per generation, and each of these
mutations has a probability $\gamma$ per generation to
disappear. While the case $\gamma=0 $ is known as (the
Fleming--Viot version of) \textit{Muller's ratchet}, the case $\gamma>
0$ is associated with \textit{compensatory mutations} in the biological
literature. We show that the system has a~unique weak solution. In
the absence of random fluctuations in type frequencies (i.e., for
the so-called infinite population limit) we obtain the solution in a
closed form by analyzing a probabilistic particle system and show
that for $\gamma>0$, the unique equilibrium state is the Poisson
distribution with parameter $\lambda/(\gamma+ \alpha)$.
\end{abstract}

% KEYWORDS
%


\end{frontmatter}

%s1 ###
\section{Introduction and outline}
We study a multitype Wright--Fisher SDE (or \textit{Fleming--Viot process})
of the form
%
%e1 ###
{\renewcommand{\theequation}{$\ast$}
\begin{eqnarray}
\label{eqSDEs}
d X_k &\!=\!& \Biggl(\! \alpha\Biggl(\!\sum_{\ell=0}^\infty(\ell-k)X_\ell
\Biggr) X_k + \lambda(X_{k-1} - X_k) +
\gamma\bigl((k+1)X_{k+1} - kX_k\bigr)\!\Biggr)\,dt
\nonumber\hspace*{-14pt}
\\[-4pt]
\\[-12pt]
\nonumber
&&{}+ \sum_{\ell\neq
k}\sqrt{\frac1N X_k X_\ell}   \,dW_{k\ell},
\end{eqnarray}
for $k=0,1,\ldots$ with $X_{-1}:=0$, $X_0, X_1,\ldots\geq0$ and
$\sum_{k=0}^\infty X_k=1$. Here $\alpha$, $\lambda$ and~$\gamma$ are
(small) nonnegative constants, $N$ is a (large) number [or equals
infinity in which case the last term on the right-hand side of (\ref{eqSDEs})
vanishes], and $(W_{k\ell})_{k>\ell}$ is a family of independent
Brownian motions with $W_{k\ell} = -W_{\ell k}$.

The interest in this system comes from population genetics (see
Section~\ref{shist} for some background). The
equations~\eqref{eqSDEs} provide a diffusion approximation of the
evolution of the \textit{type frequencies} $X_k$, $k \in\mathbb N_0$ in
a population of constant size that consists of a large number $N$ of
individuals. The type $k$ of an individual is given by the number of
deleterious mutations it carries. The \textit{fitness} (which is
proportional to the average number of offspring) of a~type-$k$
individual is proportional to $(1-\alpha)^k \approx1-\alpha k$, where
$\alpha$ is a (small) positive number called the \textit{selection
coefficient}. The parameter $\lambda$ is the expected number of
additional mutations that accumulate per individual and generation,
and for each of the mutations present, $\gamma$ is the probability
that this mutation disappears in one generation.

In this work we will not be concerned with proving the convergence of
the discrete-generation dynamics to the diffusion
approximation. Still, in Section~\ref{sssimul} we will use the
discrete generation scheme as just described in order to present a few
simulation results which illustrate how certain functionals of the
solution of~\eqref{eqSDEs} [in particular the mean and the variance
of the probability vector $(X_k)_{k=0,1,2,\ldots}$] depend on the model
parameters.

Theorem~\ref{T1} in Section \ref{ssstates} states that
\eqref{eqSDEs} has a unique weak solution. Note that \eqref{eqSDEs}
is an infinite-dimensional SDE with an \textit{unbounded nonlinear drift
coefficient}. Related existence and uniqueness results were obtained
by \citet{EthierShiga2000}. However, these authors only cover the case
of parent-independent mutation and not the situation
of~\eqref{eqSDEs}.

Theorem~\ref{T2} in Section \ref{ssNinf} gives the explicit solution
of \eqref{eqSDEs} in the (deterministic) case $N=\infty$. This
extends results from \citet{Haigh1978} and
\citet{EtheridgePfaffelhuberWakolbinger2007} for the case
$\gamma=0$. In particular, we show that the Poisson weights with
parameter $\lambda/(\gamma+ \alpha)$ constitute the only equilibrium
state of~\eqref{eqSDEs} for $N=\infty$. The proofs of
Theorems~\ref{T1} and~\ref{T2} are given in Sections~\ref{Sproof1}
and~\ref{Sproof2}, respectively. An essential step in the proof of
Theorem~\ref{T2} is Proposition~\ref{P2} which in the case $N=\infty$
provides the solution of~\eqref{eqSDEs} in terms of a probabilistic
particle system.

%s2 ###
\section{History and background of the model} \label{shist} For
$\gamma=0$, the system \eqref{eqSDEs} is known as (the Fleming--Viot
version of) Muller's ratchet, a population genetic model introduced by
Hermann \citet{Muller1964}: A clonal population of fixed size
reproduces randomly. Each individual carries a number of mutations,
all of which are assumed to be deleterious. Fitness decreases linearly
with the number of mutations. The offspring of an individual has a
small chance to gain a new deleterious mutation. In particular, any
offspring of an individual carries\vadjust{\goodbreak} at least as many mutations as the
parent, and mutation is an irreversible process. Hence, eventually the
ratchet will \textit{click} in the sense that the fittest type
will irreversibly disappear from the population. In this way, the
mutation process drives the population to a larger number of
deleterious mutations while selection acts in the opposite direction,
leading to a~form of mutation-selection
quasi-balance. \citet{GabrielLynchBuerger1993} consider a related model
of a clonally reproducing population in which the evolution of the
population size is coupled with the mean fitness of the population,
eventually leading to extinction of the population. The prediction of
this \textit{mutational meltdown} requires information on the rate at
which deleterious mutations accumulate in the
population~[\citet{Loewe2006}], that is, on the \textit{rate of Muller's
ratchet}.

Several quantitative treatments of Muller's ratchet have already been
given
[\citet{Haigh1978},
\citet{StephanChaoSmale1993},
\citet{Gessler1995},
\citet{HiggsWoodcock1995},
\citet{GordoCharlesworth2000},
\citet{MaiaBotelhoFontanari2003},
\citet{RouzineWakeleyCoffin2003},
\citet{EtheridgePfaffelhuberWakolbinger2007},
\citet{pmid18689884},
\citet{WaxmanLoewe2010},
\citet{AudiffrenPardoux2011}]. The most
interesting question concerns the rate of Muller's ratchet. This has
so far only been studied by simulations, or approximations which seem
ad hoc.

We study an extension of Muller's ratchet where deleterious mutations
are allowed to be compensated by (back-)mutations. It is important to
note that such \textit{compensatory mutations} are different from
\textit{beneficial mutations}, although both increase the fitness of an
individual. The latter are usually assumed to have an effect that does
not depend on the genetic background. In contrast, compensatory
mutations can only remove the effects of previously gained deleterious
mutations. The possibility of such compensatory mutations was
discussed already by \citet{Haigh1978} [see also \citet{MaynardSmith1978}]. He
argued that they rarely occur in
realistic parameter ranges because the deleterious mutation rate is
proportional to the full length of the genome of a clonally
reproducing individual, while the compensatory mutation rate scales
with the length of a single base within the full genome. Therefore, he
concluded that compensatory mutations are too rare to halt the
accumulation of deleterious mutations in realistic parameter
ranges. However, when several deleterious mutations are gained, the
total rate of accumulation of deleterious mutations increases and may
therefore halt the ratchet. An equilibrium is approached where a
certain number of deleterious mutations is fixed. If this number is
large enough, these may lead to extinction of the population. While
\citet{AntezanaHudson1997} argue that the effects of compensatory
mutations can be an important factor for small viruses,
\citet{Loewe2006} concludes that compensatory mutations are still too
rare to halt the mutational meltdown of human mitochondria.

Clearly, the relevance of compensatory mutations is greatest for
species with a short genome and a high mutation rate. One of the most
extreme groups in these respects are\vadjust{\goodbreak} RNA viruses (for which the genome
length is of the order of $10^3$ to $10^4$ bases and the per base
mutation rate is around $10^{-5}$ to $10^{-4}$). As discussed in
\citet{pmid2247152}, back mutations can hardly stop Muller's ratchet
even in this case. We will come back to this numerical example in
Section~\ref{sssimul} below.

The relevance of Muller's ratchet with compensatory mutations is
supported by the fact that a deleterious mutation might be compensated
not only by a back mutation that occurs at the same genomic position.
As discussed by \citet{WagnerGabriel1990}, restoring the function of a
gene which was subject to a mutation is as well possible by mutating a
second site within this gene, or even within a gene at some other
locus. \citet{Maisnier2004} give the following generalizations of
(single-base) compensatory mutations: (i)~point mutations which
restore the RNA secondary structure of a gene or the protein
structure, (ii)~an up-regulation of gene expression of the mutated
gene, (iii)~a mutation in another gene restoring the structure of a
multi-unit protein complex and (iv) a bypass mechanism where the
function of the mutated is taken over by another gene.

Various examples give clear evidence for the existence of compensatory
mutations. It has been shown by \citet{PoonChao2005} that a deleterious
mutation in the DNA bacteriophage phiX174 can be compensated by about
nine different intragenic compensatory mutations. This implies that
the rate of compensatory mutations can be high enough to halt
accumulation of deleterious mutations under realistic scenarios. In
fact, compensatory mutations have been observed in various
species. \citet{Howe2008} showed that deletions in protein-coding
regions of the mitochondrial genome in \textit{Caenorhabditis briggsae}
lead to heteroplasmy, a~severe factor in mitochondrial diseases. They
also found compensatory mutations leading to a~decrease in
heteroplasmy. Mutations for antibiotic resistance of bacteria are
known to be deleterious in a wild-type population. Fitness can be
increased by a compensatory mutation [see, e.g., \citet{Handel2006}].
Plastid genomes of mosses are studied in
\citet{Maier2008BMC-Biol18755031}. Here, it is suggested that
deleterious mutations may be compensated by RNA editing, a mechanism
by which the base $C$ in DNA is transcribed to $U$ on the level of RNA
for specific bases in the genome.

All these examples indicate that the role of compensatory mutations
should be taken into account. A relevant question to be addressed in
future research is which parameter constellations (of the selection
coefficient, the mutation rate, the compensatory mutation rate and the
population size) can halt the ratchet before the mutational meltdown
leads to extinction of the population.

%s3 ###
\section{Results}
\label{sresults}
We show that for finite $N$ the system~\eqref{eqSDEs} has a unique
weak solution (Theorem~\ref{T1}). For the system~\eqref{eqSDEs}
without noise (i.e., the case $N=\infty$) we provide in Theorem \ref{T2} the
explicit form of the solution as well as the equilibrium state.\vadjust{\goodbreak} For
this we use a stochastic particle model (including accumulation and
loss of mutations, as well as a state-dependent death rate of the
particles) and show in Proposition~\ref{P2} that a solution
of~\eqref{eqSDEs} with $N=\infty$ is given by the distribution of the
particle system conditioned on nonextinction. After stating the
theorems, we compare in Section~\ref{sssimul} the cases of large $N$
with the theoretical considerations for $N=\infty$ using simulations.

%s3.1 ###
\subsection{Existence and uniqueness}
\label{ssstates}
The system \eqref{eqSDEs} of Muller's ratchet with compensatory
mutations takes values in the space of probability vectors indexed by~$\mathbb N_0$,
that is, sequences whose entries are probability weights
on~$\mathbb N_0$. We restrict the state space to the subset of
probability vectors with finite exponential moment of a certain order,
and show uniqueness in this space. Throughout, we abbreviate
$\underline x := (x_0, x_1,\ldots) \in\mathbb R_+^{\mathbb N_0}$.

\begin{definition}[Simplex]
The \label{simp} \textit{infinite-dimensional simplex} is given by
\renewcommand{\theequation}{\arabic{section}.\arabic{equation}}
\begin{equation}
\mathbb S:=\Biggl\{\underline x\in\mathbb R_+^{\mathbb N_0}\dvtx
\sum_{k=0}^\infty x_k=1\Biggr\}.
\end{equation}
Moreover, for $\xi>0$, set
\begin{equation}\label{eqhdef}
h_\xi(\underline x) := \sum_{k=0}^\infty x_k e^{\xi k}
\end{equation}
and consider elements of $\mathbb S$ with $\xi$th exponential
moment, forming the space
\begin{equation}\label{eqdefScirc}
{\mathbb S}_\xi:= \{ \underline x \in\mathbb S\dvtx
h_\xi(\underline x) < \infty\}.
\end{equation}
\end{definition}

\begin{remark}[(Topology on $\mathbb S_\xi$)]
We note that \label{remtopo}
\renewcommand{\theequation}{\arabic{section}.\arabic{equation}}
\begin{equation}
\label{eqrmetrix}
r(\underline x, \underline y) := \sum_{k=0}^\infty e^{\xi
k}|x_k-y_k|, \qquad \underline x, \underline y \in\mathbb S_\xi,
\end{equation}
defines a complete and separable metric on $ \mathbb S_\xi$.
\end{remark}

\begin{theorem}[(Well-posedness of Fleming--Viot system)] \label{T1}
Let \ $\underline x \in\mathbb S_\xi$ for some $\xi>0$.
Then, for $N\in(0,\infty)$, $\alpha,\lambda,\gamma\in[0,\infty)$,
the system~\eqref{eqSDEs} starting in $\underline X(0)=\underline
x$ has a unique $\mathbb S$-valued weak solution $\mathcal X =
(\underline X(t))_{t\geq0}$, taking values in the space $\mathcal
C_{\mathbb S_\xi}([0,\infty))$ of continuous functions on $\mathbb
S_\xi$.
\end{theorem}

In the sequel, we will refer to the process $\mathcal X$ as
\textit{Muller's ratchet with compensatory mutations} with selection
coefficient $\alpha$, mutation rate $\lambda$, compensatory mutation
rate $\gamma$ and population size $N$.

\begin{remark}[(Population size $N$)]
Resampling \label{Popsize} models are usually studied either for a
finite population of constant size $N$ (e.g., using a Wright--Fisher
model), or in the large population limit with a suitable rescaling
of time, leading to Fleming--Viot processes. For a bounded fitness\vadjust{\goodbreak}
function and a compact type space, it is well known that a sequence
of (discrete time) Wright--Fisher processes, indexed by $N$,
converges weakly to a Fleming--Viot process (or Wright--Fisher
diffusion) if the selection and mutation coefficients are scaled
down by $N$ and one unit of time is taken as $N$ generations; see,
for example, \citet{EthierKurtz1993}. In our situation it may thus be
expected (though we do not prove this claim here) that for large~$N$
and for $\alpha N$, $\lambda N$ and $\gamma N$ of order one, the
Wright--Fisher process described in Section~\ref{sssimul}, run with
a time unit of $N$ generations, is close to the solution
of~\eqref{eqSDEs}, with~$1$ instead of $\sqrt{1/N}$ as the
coefficient of the noise, and $\alpha N$, $\lambda N$ and $\gamma N$
in place of $\alpha$, $\lambda$ and $\gamma.$ However, this system
is \eqref{eqSDEs} with time speeded up by a~factor~$N$. In other
words, for large $N$, and $\alpha N$, $\lambda N$ and $\gamma N$ of
order one, the solution of \eqref{eqSDEs} should be close to the
corresponding Wright--Fisher model as introduced in
Section~\ref{sssimul}, with time unit being one generation. This is
the reason why we refer to the model parameter $N$ in
\eqref{eqSDEs} as the population size. We use this terminology in
interpreting the simulation results for the Wright--Fisher model in
Section~\ref{sssimul}.
\end{remark}

\begin{remark}[(Connection to previous work for $\gamma=0$)]
For the case $\mu= 0$, variants of Theorem \ref{T1} appear in
\citet{Cuthbertson2007} and in \citet{AudiffrenPardoux2011}. The
latter makes (in the terminology of our Theorem~\ref{T1}) the
special choice $\xi= \alpha N$ and refers to \citet{Audiffren2011}
for the proof. \citet{Cuthbertson2007} treats also the case of
$\alpha< 0$, assuming the existence of all exponential moments of
the initial state.
\end{remark}

\begin{remark}[(Strategy of the proof of Theorem \protect\ref{T1})]
For $\alpha=0$, it follows from classical theory [\citet{Dawson1993}, Theorem
5.4.1] that \eqref{eqSDEs} has a unique weak
solution. The same is true if the selection term $\alpha(\sum_{\ell
=0}^\infty(\ell-k)X_\ell)X_k$ is replaced by a bounded function
of $\underline X$. This can be shown by a Cameron--Martin--Girsanov
change of measure from the case $\alpha=0$, using similar arguments
as in \citet{EthierShiga2000}. So, the main difficulty in the proof is
to deal with the unbounded selection term. This is overcome by
showing that the change of measure still works when using $\mathbb
S_\xi$ as the state space for~$\mathcal X$.
\end{remark}

\begin{remark}[{[Strong solution of \eqref{eqSDEs}]}]
Theorem~\ref{T1} gives existence and uniqueness of \textit{weak
solutions} of~\eqref{eqSDEs}. To the best of our knowledge, a
result on uniqueness of \textit{strong solutions} so far is not
available even in the case $\gamma=\lambda=\alpha=0$. The reason why
general theory does not apply in this multidimensional situation is
that the diffusion term $\sqrt{X_kX_\ell}$ is only H\"older rather
than Lipschitz continuous. However, let us mention two related
results:
\begin{longlist}[(ii)]
\item[(i)] \citet{ShigaShimizu1980} provide, in their Theorem 3.4,
existence and uniqueness of strong solutions for a class of SDE's
which are similar to our system~\eqref{eqSDEs}. One may
conjecture that this theorem is also valid for the drift term\vadjust{\goodbreak}
appearing in our system \eqref{eqSDEs}. This would then give an
alternative proof of our Theorem~\ref{T1}. However, the diffusion
term in the SDE considered in \citet{ShigaShimizu1980} is assumed
to have a lower triangular form, which seems to be a tribute to
the mathematical technique rather than to a biological
interpretability from a ``strong'' (i.e., realization-wise) point
of view.

\item[(ii)] Recently, \citet{DawsonLi2010} [see their equation~(4.1)]
studied strong existence and uniqueness for a related system of
stochastic flows. Here, white noise on $[0,\infty) \times[0,1]$
is used to model the reproduction of the individuals in the
(unstructured) population, irrespective of their type.
\end{longlist}
\end{remark}

%s3.2 ###
\subsection{\texorpdfstring{The case $N=\infty$}{The case N=infinity}}
\label{ssNinf}
This case (which is not included in Theorem~\ref{T1}) leads to a
deterministic dynamics. For $\gamma=0$, \citet{Haigh1978} was the first
to obtain results on the deterministic evolution of $\mathcal X$ in a
discrete time setting. These results were later refined by
\citet{MaiaBotelhoFontanari2003}. Here, we work with continuous time,
and our next theorem generalizes Proposition~4.1 in
\citet{EtheridgePfaffelhuberWakolbinger2007} to the case $\gamma>0$. We
are dealing with the system
\renewcommand{\theequation}{\arabic{section}.\arabic{equation}}
\begin{equation}\label{eqode1}\qquad
\dot x_k = \alpha\Biggl(\sum_{\ell=0}^\infty(\ell-k)x_\ell\Biggr) x_k +
\lambda(x_{k-1} - x_k)+ \gamma\bigl((k+1)x_{k+1} - kx_k\bigr)
\end{equation}
for $k=0,1,2,\ldots$ with $x_{-1}:=0$ and $\sum_{k=0}^\infty x_k=1$.

\begin{theorem}
\label{T2}
Let $\alpha, \lambda, \gamma\in[0,\infty)$ and $\underline x(0)
\in\mathbb S_\xi$ for some $\xi>0$. Then
system~\eqref{eqode1} has a unique $\mathbb S$-valued solution
$(\underline x(t))_{t\geq0}$ which takes values in~$\mathbb S_\xi$.
It is given by
\renewcommand{\theequation}{\arabic{section}.\arabic{equation}}
\begin{eqnarray}\label{eqT11}
x_k(t) & =& \Biggl( \sum _{i=0}^\infty x_i(0)
\sum _{j=0}^{i\wedge k} \pmatrix{{i}\cr{j}} \bigl({\bigl(\gamma\bigl(
1-e^{-(\alpha+\gamma)t}\bigr)\bigr)}/{(\alpha+\gamma)}\bigr)^{i-j}
\nonumber\\
&&\hspace*{40pt}\qquad{}\times e^{-j(\alpha+\gamma)t} \bigl({1}/{(k-j)!}\bigr)\nonumber\\
&&\hspace*{71pt}{}\times\bigl({\bigl(\lambda\bigl(1-e^{-(\alpha+\gamma)t}\bigr)\bigr)}/{(\alpha+\gamma)}
\bigr)^{k-j}\Biggr)\\
&&{}\Bigg/\Biggl( \sum _{i=0}^\infty x_i(0)
\bigl({\gamma}/{(\alpha+\gamma)} - {\alpha}/{(\alpha+\gamma)}
e^{-(\alpha+\gamma)t}\bigr)^i\nonumber \\
&&\hspace*{44pt}\qquad{}\times
\exp\bigl({\lambda}/{(\alpha+\gamma)}\bigl(1-e^{-(\alpha+\gamma)t}\bigr)
\bigr)\Biggr).\nonumber
\end{eqnarray}
In particular, if either $\gamma>0$ or $x_0(0)>0$, then
\begin{equation}\label{eqT12}
x_k(t)  \rightarrow{t\to\infty}
\frac{e^{-\lambda/(\alpha+\gamma)}}{k!}
\cdot\biggl(\frac{\lambda}{\alpha+\gamma}\biggr)^k, \qquad  k = 0,1,2,
\ldots;
\end{equation}
that is, the limiting state as $t \to\infty$ is the vector of Poisson
weights with
parameter $\lambda/(\alpha+\gamma)$.
\end{theorem}

\begin{remark}[(Equilibria)]
In the case $\gamma= 0$ it is known already from the work of
\citet{Haigh1978} that the vector of Poisson weights with parameter
$\lambda/\alpha$ is an equilibrium state. Moreover, Poisson states
which are shifted by $k = 1, 2, \ldots$ are equilibria as well. This
is in contrast to the case $\gamma> 0$ where only a single
equilibrium state exists. Moreover, this equilibrium state depends
on the model parameters only through the value of
$\lambda/(\alpha+\gamma)$. It is worth noting, however, that for
finite $N$ the distribution of the process does not merely depend on
$\lambda/(\alpha+\gamma)$. See Figure~\ref{fig3} for a simulation
study of this feature for $N < \infty$.

Here is a heuristic argument why the equilibrium state is given by a
Poisson distribution with parameter $\lambda/(\gamma+ \alpha)$ in
the case $N=\infty$. Consider the number of mutations along a single
line which are accumulated at rate $\lambda$, each of which is
compensated at rate $\gamma$, and if the line carries $k$ mutations,
the line is killed at rate~$\alpha k$. We will see in
Proposition~\ref{P2} that the equilibrium distribution
of~\eqref{eqode1} equals the quasi-equilibrium for this Markov
chain, that is, the equilibrium distribution of mutations given the
line has not been killed. In the case $\alpha= 0$ (and $\lambda,
\gamma> 0$), the number of mutations can be seen as an
$M/M/\infty$-queueing system: In equilibrium an individual carries
Pois($ \lambda/ \gamma$)-many mutations. These mutations have
i.i.d.\ Exp($\gamma$)-distributed ages. In the case $\alpha> 0$,
the equilibrium distribution can be constructed from that for
$\alpha=0$ by conditioning on the event that none of the present
lines experienced the killing caused by the mutations it carries.
Since each mutation has an independent Exp$(\gamma)$ distributed
age, this gives the chance $(\frac
{\gamma}{\gamma+\alpha})^k$ for survival of the line. The
claim then results from the following elementary fact: Consider a
population with Poisson($\beta$)-distributed size with $\beta=
\lambda/\gamma$. Then, conditional under the event that this
population remains completely intact under a thinning with survival
probability $p=\gamma/(\alpha+\gamma)$, the size of the population
is Poisson($\beta p$) distributed.
\end{remark}

\begin{remark}[(Connection to the rate of adaptation)]
Although the strategy of our proof requires that $\alpha\geq0$
(i.e., the mutations are deleterious), it can be shown by taking the
time-derivative of the right-hand side of~\eqref{eqT11} that this
equation is a solution for $\alpha<0$ as well. This model is
frequently termed \textit{rate of adaptation} and has gained some
interest in the case $\gamma=0$ and $N<\infty$
[\citet{GerrishLenski98},
\citet{DesaiFisher2007},
\citet{ParkKrug2007},
\citet{YuEtheridgeCuthbertson2010}].

Taking $\alpha<0$ in our model, all mutations are beneficial, and
$\gamma$ is the rate by which any beneficial mutation is
compensated. Interestingly, only in the case \mbox{$|\alpha| < \gamma$}
(i.e., selection is weaker than the compensatory mutation rate) an
equilibrium state exists, and is still Poisson with parameter
$\lambda/(\gamma- |\alpha|)$. In the case $|\alpha| \geq\gamma$,
no equilibrium exists because new beneficial mutations spread
through the population quicker than compensatory mutations can halt
this process. It will be interesting to investigate the switch between these
two scenarios in the case of finite $N$.
\end{remark}

%s3.3 ###
\subsection{Simulations}
\label{sssimul}
We use simulations based on a discrete Wright--Fisher model to study
the evolution of the mean fitness, and to investigate the dependence
of the mean and the variance of the type frequency distribution on the
model parameters. Fixing a population size $N$, this model is a
discrete time Markov chain $(\underline Y(t))_{t=0,1,2,\ldots}$ taking
values in $\{\underline y \in\mathbb S\dvtx  N\underline y \in\mathbb
N_0^{\mathbb N_0}\}$ and such that
\[
\mathbf P\bigl(\underline Y(t+1) = \underline y|\underline Y(t)\bigr) = \pmatrix{{N}\vspace*{2pt}\cr{Ny_0
  Ny_1   Ny_2 \cdots} }\prod_{j=0}^\infty p_j^{Ny_j},
\]
where
\begin{eqnarray*}
\mathrm{(i)}&&\quad  \widetilde p_j  = \frac{(1-\alpha)^j
Y_j(t)}{\sum_{k=0}^\infty
(1-\alpha)^k Y_k(t)},\\
\mathrm{(ii)}&&\quad  \widehat p_j  = \sum_{m=j}^\infty\widetilde p_m
\pmatrix{{m}\vspace*{2pt}\cr{j}}\gamma^{m-j}
(1-\gamma)^{j},\\
\mathrm{(iii)}&&\quad  p_j  = \sum_{l=0}^j \widehat p_l
e^{-\lambda}\frac{\lambda^{j-l}}{(j-l)!}
\end{eqnarray*}
for \textit{small} parameters $\alpha, \lambda$ and $\gamma$.
The sampling weights $(p_j)_{j=0,1,\ldots}$ describe selection,
mutation and compensatory mutation. The idea in this scheme (which is
standard in population genetics) is that: (i) any individual produces a
large number of gametes, but an individual with $k$ deleterious
mutations only contributes a number proportional to $(1-\alpha)^k$ to
the gamete pool; (ii) every deleterious mutation has a small,
independent chance $\gamma$ to be removed while the gamete is built;
(iii)~the number of new deleterious mutations is Poisson
distributed with parameter~$\lambda$. After building these gametes,
$N$ individuals are randomly chosen from the gamete pool to form the
next generation. Since $\alpha, \gamma$ and $\lambda$ are assumed to
be small, the order in which the three mechanisms (i), (ii), (iii)
come into play is negligible. (E.g., if we would assume---in contrast
to our simulation scheme above---that compensatory mutations arise
before gametes are built proportionally to the relative fitness of
individuals. Then an individual with a high number of deleterious
mutations would produce slightly more gametes than in our simulation
scheme.) For our simulations, the working hypothesis is that
$(\underline Y(Nt))_{t\geq0}$ behaves similarly to $\mathcal X
=(\underline X(Nt))_{t\geq0}$ where $\mathcal X$ is the solution
of~\eqref{eqSDEs} with parameter $N$; see Remark~\ref{Popsize}.

%f1 ###
\begin{figure}

\includegraphics{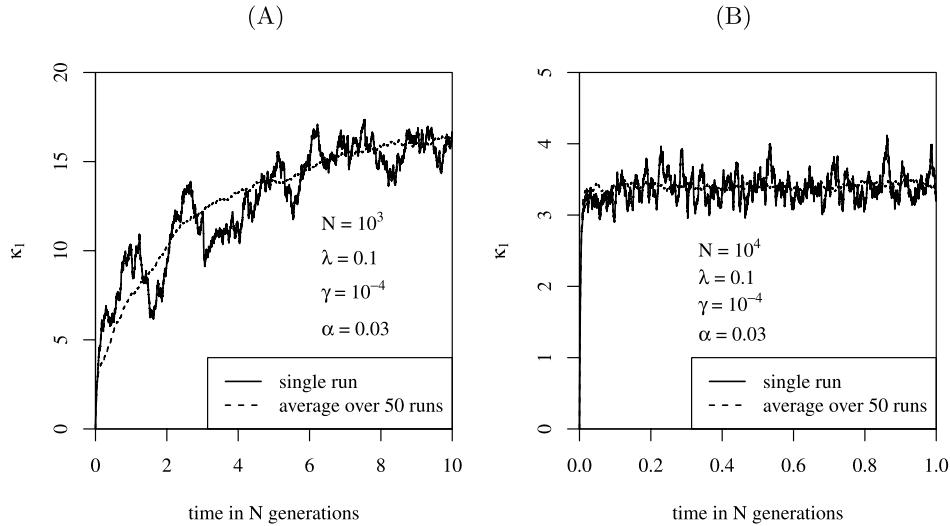}

\caption{The evolution of the average number of
deleterious mutations $\kappa_1$ is plotted. In addition, the
single path is compared to the average over 50 different
simulations. \textup{(A)}~A~parameter combination where Muller's ratchet
without compensatory mutations (i.e., $\gamma=0$) clicks
frequently, while it clicks much less frequent in
\textup{(B)}.}\label{fig1}\vspace*{3pt}
\end{figure}

We simulated $(\underline Y(Nt))_{t\geq0}$ for various combinations
of $N, \alpha, \lambda$ and $\gamma$, starting with $\underline
Y(0)=\delta_0$; that is, no deleterious mutations are present at
start. Since in reality compensatory mutations are less probable than
mutations, we mostly simulate scenarios with $\gamma\ll\lambda$. (For
the biological background of this assumption, see
Section~\ref{shist}.) Hence, our simulations can be considered as a
small perturbation of the case $\gamma=0$, the case of Muller's
ratchet (without compensatory mutations). We compare scenarios where
Muller's ratchet clicks rarely with others in which it clicks more
frequently. For example, in Figure~\ref{fig1}(A) we use $N=10^3,
\lambda=0.1, \alpha=0.03$ where the ratchet has about 5.7 clicks in
$N$ generations. In Figure~\ref{fig1}(B) we use $N=10^4$ where the
ratchet has only about 0.34 clicks in $N$ generations. Both figures
show the initial phase of the simulation for a small compensatory
mutation rate of $\gamma=10^{-4}$. Recall that Theorem~\ref{T2}
predicts an equilibrium number of $\lambda/(\alpha+\gamma) \approx3.3$
deleterious mutations in the case $N=\infty$. This value is reflected
in our simulations only in Figure \ref{fig1}(B) where Muller's ratchet clicks
rarely. In Figure~\ref{fig1}(A), not only is the average number of deleterious
mutations much larger than the prediction from Theorem~\ref{T2}, but
also the fluctuations are much larger than in Figure~\ref{fig1}(B). However, in
both parameter constellations we see that the accumulation of
deleterious mutations by Muller's ratchet is slowed down (and sooner
or later halted) due to the compensatory mutations.

%f2 ###
\begin{figure}

\includegraphics{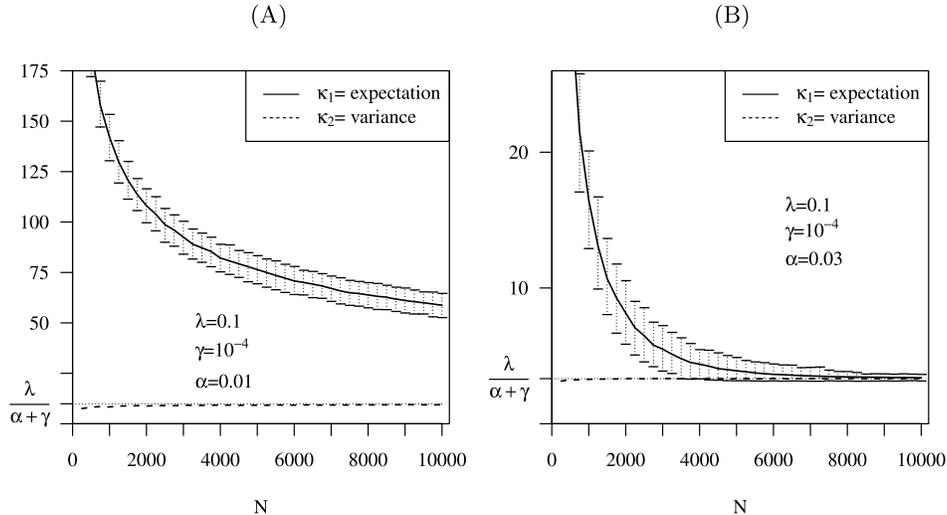}

\caption{The empirical distribution of $\kappa_1$ and
$\kappa_2$ are evaluated between generations $5\cdot10^2 N$ and
$10^3 N$. The plot for $\kappa_1$ includes the resulting 10\% and
90\% quantiles. In absence of compensatory mutations and with
$N=10^4$, the same parameters lead to approximately 152 clicks per
$N$ time units for \textup{(A)}, while 0.34 clicks per $N$ time
units are obtained for~\textup{(B)}.}\label{fig2}
\end{figure}

Figure~\ref{fig2} illustrates for a finite $N$, how far the mean and
variance of the number of deleterious mutations deviate from those
Poisson distribution, which appears in Theorem~\ref{T2} for the case
$N=\infty$. Again, we see that for fixed $\alpha, \lambda$ and small
compensatory mutation rate $\gamma$, the equilibrium for $\kappa_1$ is
close to $\lambda/(\alpha+\gamma)$ only if the ratchet without
compensatory mutations ($\gamma=0$) does not click too often. If
$N=10^4$ in Figure~\ref{fig2}(A), there are approximately 152 clicks
in $N$ generations in the absence of compensatory mutations, while in
Figure~\ref{fig2}(B), this rate is much lower, approximately 0.34
clicks per $N$ generations [using the same parameter values $\alpha$,
$\lambda$ and $\gamma$ as in Figure~\ref{fig1}(B)]. These examples
show that compensatory mutations halt the ratchet quite efficiently.
Note that the parameter values $\lambda=0.1$ and $\gamma=10^{-4}$ fit
to the evolution of RNA viruses, for example, for a genome of length $10^3$
bases, if the per base mutation rate is $10^{-4}$ and a population of
size $10^4$. As our simulations show, the ratchet is halted, provided
the selection coefficient is large enough. This is in some contrast to
\citet{pmid2247152} who argues that compensatory mutations are too rare
in RNA viruses to halt the ratchet.

Another surprising fact in Figure~\ref{fig2}(A) is that the empirical variance
of the number of deleterious mutations in the population is always
close to the prediction of $\lambda/(\alpha+\gamma)$ from the Poisson
state appearing in Theorem~\ref{T2}. This would be compatible
with the hypothesis that the type frequencies for the Wright--Fisher
model are (in equilibrium) close to a shifted Poisson
distribution. The detailed study of the amount of this shift, in
particular, for a~parameter constellation for which $\gamma=0$ leads to
frequent clicks, is a~delicate issue. Its answer certainly depends on
the rate of clicking of Muller's ratchet without compensatory
mutations, a problem which remains unsolved until now.

Yet another interesting feature seen in Theorem~\ref{T2} is the
symmetric dependence on $\alpha$ and $\gamma$ of the equilibrium
state. We checked in which direction this symmetry is violated for
finite $N$. As seen from Figure~\ref{fig3}, compensatory mutations can
halt the ratchet more efficiently than selection. The reason is that
compensatory mutations reduce the number of mutations no matter how
many mutations are fixed in the population, whereas the number of
fixed mutations cannot decrease due to selection.

%f3 ###
\begin{figure}

\includegraphics{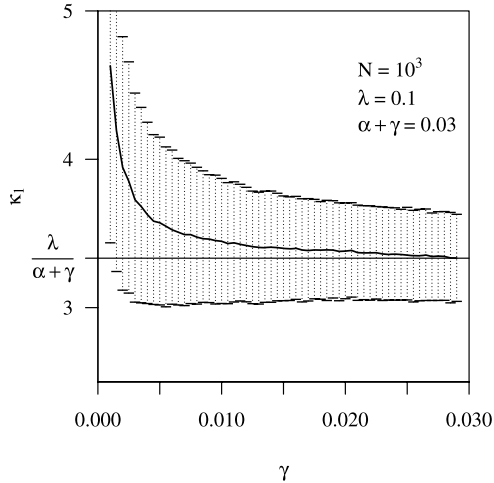}

\caption{The parameter $\lambda/(\alpha+\gamma)$ of the
Poisson equilibrium state in the case $N=\infty$ is
symmetric in $\alpha$ and $\gamma$ (see Theorem~\protect\ref{T2}). We fix
$\alpha+\gamma$ and see that the average number of deleterious
mutations is higher for low values of $\gamma$. Again, the 10\%
and 90\% quantiles are given.}\label{fig3}
\end{figure}

%s4 ###
\section{\texorpdfstring{Proof of Theorem~\protect\ref{T1}}{Proof of Theorem 1}}
\label{Sproof1}
Our approach is inspired by \citet{EthierShiga2000}, who deal with the
case of unbounded selection if mutation is parent-independent. In
order to prove that~\eqref{eqSDEs} has a unique weak solution, we use
the equivalent formulation by a martingale problem and show its
well-posedness (see Proposition~\ref{P34}). We provide bounds on
exponential moments for any solution of the martingale problem
associated with the generator of~\eqref{eqSDEs} in
Lemma~\ref{lexpBounds}. The central step is Proposition~\ref{P2}
which provides a~Girsanov change of measure by which solutions of the
martingale problem for $\alpha=0$ are transformed to solutions for any
$\alpha>0$. Proposition~\ref{P34} and Theorem~\ref{T1} then follow
because the martingale problem for $\alpha=0$ is well posed, and can
be transformed to a solution for the martingale problem for $\alpha
>0$ which also solves \eqref{eqSDEs}. This shows
existence. Uniqueness again follows by using a Girsanov transform.

%s4.1 ###
\subsection{Martingale problem}
\label{ssmp}
We start by defining the generator for the Fleming--Viot system of
Muller's ratchet with compensatory mutations. The unboundedness in
the selection term of this generator requires particular care in the
analysis of the corresponding martingale problem. First we fix some
notation.

\begin{remark}[(Notation)]
For a complete and separable metric space $(\mathbb E,r)$, we denote
by $\mathcal P(\mathbb E)$ the space of probability measures on (the
Borel sets of) $\mathbb E$, and by $\mathcal M(\mathbb E)$
[resp., $\mathcal B(\mathbb E)$] the space of real-valued, measurable
(and bounded) functions. If $\mathbb E\subseteq\mathbb R^{\mathbb
N_0}$, we let $\mathcal C^k(\mathbb E)$ ($\mathcal C_b^k(\mathbb
E)$) be the (bounded), $k$ times partially continuously
differentiable functions (with bounded derivatives). Partial
derivatives of $f\in\mathcal C^2(\mathbb E)$, $\mathbb E\subseteq
\mathbb R^{\mathbb N_0}$, will be denoted by
\renewcommand{\theequation}{\arabic{section}.\arabic{equation}}
\begin{equation}\label{eqdefPartDer}
f_k := \frac{\partial f}{\partial x_k}, \qquad  f_{k\ell} :=
\frac{\partial f^2}{\partial x_k\, \partial x_\ell},\qquad
k,\ell=0,1,2,\ldots.
\end{equation}
\end{remark}

\begin{definition}[(Martingale problem)]
Let \label{defmp} $(\mathbb E,r)$ be a complete and separable
metric space, $x\in\mathbb E$, $\mathcal F \subseteq{\mathcal
M}(\mathbb E)$ and $G$ a linear operator on ${\mathcal M}(\mathbb
E)$ with domain~$\mathcal F$. A (distribution $\mathbf P$ of an)
$\mathbb E$-valued stochastic process $\mathcal X=(X_t)_{t\geq0}$
is called a solution of the $(\mathbb E, x, G,\mathcal
F)$-martingale problem if $X_0=x$, and $\mathcal X$ has paths in the
space ${\mathcal D}_{\mathbb E}([0,\infty))$, almost surely, and for
all $f\in\mathcal F$,
%
%e2 ###
\renewcommand{\theequation}{\arabic{section}.\arabic{equation}}
\begin{equation}\label{13def}
\biggl(f(X_t)-f(X_0) - \int_0^t G f(X_s) \,ds \biggr)_{t\geq0}
\end{equation}
is a $\mathbf{P}$-martingale with respect to the canonical
filtration. Moreover, the $(\mathbb E, x, G,  \mathcal F)$-martingale
problem is said to be well posed if there is a unique solution
$\mathbf{P}.$
\end{definition}

For a fixed $\xi> 0$, our state space will be ($\mathbb S_\xi, r$);
cf. Definition \ref{simp} and Remark~\ref{remtopo}. We now specify
the generator and its domain.

% % \newpage

\begin{definition}[(Generator for Fleming--Viot
system)]\label{defFVMR}
(1) On $\mathbb S$, consider functions of the form
%
%e3 ###
\renewcommand{\theequation}{\arabic{section}.\arabic{equation}}
\begin{eqnarray}
\label{eq122}
f(\underline x) & :=& f_{\varphi_1,\ldots,\varphi_n}(\underline x)
:= \langle\underline x, \varphi_1\rangle
\cdots\langle\underline x, \varphi_n\rangle,
\nonumber
\\[-8pt]
\\[-8pt]
\nonumber
\langle\underline x, \varphi\rangle& :=& \sum_{k=0}^\infty
x_k \varphi(k)
\end{eqnarray}
for $n=1,2,\ldots$ and $\varphi, \varphi_1,\ldots,\varphi_n \in\mathcal
M(\mathbb N_0)$. Let
%
%e4 ###
\begin{eqnarray}\label{eq121b}%
\mathcal F &:=& \mbox{the algebra generated by}\nonumber\hspace*{-35pt}
\\[-4pt]
\\[-12pt]
&&{}\{f_{\varphi_1,\ldots,\varphi_n}\dvtx  \varphi_i\in\mathcal
M(\mathbb N_0) \mbox{ with bounded support}, i=1,\ldots,n, \ n\in\mathbb
N\}.\nonumber\hspace*{-35pt}\vadjust{\goodbreak}
\end{eqnarray}

\begin{longlist}[(2)]
\item[(2)] We define the operator $G_{\mathcal X}^\alpha$ as the linear
extension of
%
%e5 ###
\begin{eqnarray}
\label{eqgen1}
G_{\mathcal X}^{\alpha} f(\underline x) & =&
G_{\mathrm{sel}}^{\alpha} f(\underline x) + G_{\mathrm{mut}}
f(\underline x) + G_{\mathrm{cmut}}
f(\underline x) + G^N_{\mathrm{res}} f(\underline x),\nonumber\\
G_{\mathrm{sel}}^{\alpha} f(\underline x) & =& \alpha
\sum_{k=0}^\infty\sum_{\ell=0}^\infty(\ell-k)x_\ell x_k
f_k(\underline x),\nonumber\\
G_{\mathrm{mut}} f(\underline x) & =& \lambda\sum_{k=0}^\infty(
x_{k-1} - x_k) f_k(\underline
x),\\
G_{\mathrm{cmut}} f(\underline x) & =& \gamma\sum_{k=0}^\infty\bigl(
(k+1)x_{k+1} - kx_k\bigr) f_k(\underline
x),\nonumber\\
G^N_{\mathrm{res}} f(\underline x) &=& \frac1{2N} \sum_{k,\ell=
0}^\infty x_k ( \delta_{k\ell} - x_\ell) f_{k\ell}(\underline
x)\nonumber
\end{eqnarray}
with $\alpha, \lambda, \gamma\in[0,\infty)$, $N\in(0,\infty)$,
for $f$ of the form~\eqref{eq122} whenever the right-hand sides
of \eqref{eqgen1} exist (which is certainly the case if
$\underline x$ has a first moment and the
$\varphi_1,\ldots,\varphi_n$ have bounded support). In particular,
for all $f \in\mathcal F$ and $\xi> 0$, the function
$G_{\mathcal X}^{\alpha} f$ is defined on $\mathbb S_\xi$.
\item[(3)] For $f = f_{\varphi_1, \ldots, \varphi_n}$, we define
$\mathcal
N_f = (N_f(t))_{t\geq0}$ by
\begin{equation}
\label{eqNf}
N_f(t) := f(\underline X(t)) - \int_0^t G_{\mathcal X}^\alpha
f(\underline X(s))\,ds
\end{equation}
whenever $G_{\mathcal X}^\alpha f(\underline X(t))$ exists for all
$t\geq0$.
\end{longlist}
\end{definition}

\begin{proposition}[(Martingale problem is well-posed in $\mathbb
S_\xi$)]
Let $\underline x\in\mathbb S_\xi$ for some $\xi>0$,
\label{P34}
$G_{\mathcal X}^\alpha$ as in~\eqref{eqgen1} and $\alpha, \lambda,
\gamma\in[0,\infty), N\in(0,\infty)$ and ${\mathcal F}$ be as
in~\eqref{eq121b}. Then the $(\mathbb S, {\underline x},
G_{\mathcal X}^\alpha, {\mathcal F})$-martingale problem is
well posed and is a~process with paths in $\mathcal C_{\mathbb
S_\xi}([0,\infty))$.
\end{proposition}

Proposition~\ref{P34} is a crucial step in the proof of
Theorem~\ref{T1}. Both proofs are carried out in
Section~\ref{SproofT1}. Now we start with bounds on exponential
moments, which will be fundamental in further proofs.

\begin{lemma}[(Bounds on exponential moments)]\label{lexpBounds}
Let $\underline x \in\mathbb S_\xi$ for some \mbox{$\xi>0$} and $\mathcal
X = (X(t))_{t\geq0}$ be a solution of the $(\mathbb S, {\underline
x},G_{\mathcal X}^\alpha, \mathcal F)$-martingale problem. Then
\renewcommand{\theequation}{\arabic{section}.\arabic{equation}}
\begin{equation}\label{eqexpBounds0}
\mathbf E[h_{\xi}(\underline X(t))] \leq h_{\xi}(\underline x)
\cdot\exp\bigl( \lambda t(e^\xi- 1)\bigr)
\end{equation}
and for all $T>0$ and $\varepsilon>0$, there is $C>0$, depending on
$T, \varepsilon, \xi$ and $\lambda$ (but not on $\alpha, \gamma, N$)
with
\begin{equation}\label{eqexpBounds1}
\mathbf P\Bigl[  \sup_{0\leq t\leq T} h_\xi(\underline
X(t))>C\Bigr] \leq\varepsilon\cdot h_\xi(\underline x).
\end{equation}
\end{lemma}

\begin{pf}
Define for $m=0,1,2,\ldots$ the function $h_{\xi,m} \in\mathcal F$ by
\[
h_{\xi,m}(\underline x) := \sum_{k=0}^m
x_k e^{\xi k} + e^{\xi m} \Biggl(1 - \sum_{k=0}^m x_k\Biggr) = e^{\xi m}
+ \sum_{k=0}^m x_k(e^{\xi k} -e^{\xi m}),
\]
and note that
\[
h_{\xi,m}(\underline x) = \sum_{k=0}^\infty
x_k e^{\xi(k\wedge m)} \qquad \mbox{for }\underline x\in\mathbb
S.
\]
First, we compute
\begin{eqnarray*}
G_{\mathrm{mut}} h_{\xi,m}(\underline x) & =& \lambda\sum_{k=0}^m
(x_{k-1} - x_k) (e^{\xi k}-e^{\xi m}) = \lambda\sum_{k=0}^{m-1}
x_k \bigl( e^{\xi(k+1)} -
e^{\xi k}\bigr) \\[-2pt]
& =& \lambda(e^\xi-1)\sum_{k=0}^{m-1}x_ke^{\xi k}\geq
0,\\[-2pt]
G_{\mathrm{cmut}} h_{\xi,m}(\underline x) & = &\gamma\sum_{k=0}^m
\bigl((k+1)x_{k+1} - kx_k\bigr) (e^{\xi k}-e^{\xi m}) \\[-2pt]
& =& \gamma
\sum_{k=1}^m kx_k \bigl( e^{\xi(k-1)} - e^{\xi k}\bigr) \leq0,
\\[-2pt]
G^\alpha_{\mathrm{sel}} h_{\xi,m}(\underline x) & =& \alpha
\sum_{k=0}^m \sum_{\ell=0}^\infty(\ell-k)x_\ell x_k (e^{\xi
k}-e^{\xi m})\\[-2pt]
& =& \alpha\sum_{k=0}^\infty
\sum_{\ell=0}^\infty(\ell-k)x_\ell x_k \bigl(e^{\xi(k\wedge
m)}-e^{\xi m}\bigr) \\[-2pt]
& =& \alpha\sum_{k=0}^\infty
\sum_{\ell=0}^\infty(\ell-k)x_\ell x_k e^{\xi(k\wedge m)} \leq
0,
\end{eqnarray*}
where the calculation for the term $G^\alpha_{\mathrm{sel}}$ holds for
$\underline x\in\mathbb S$. (For the last inequality,
assume that $Z$ is an $\mathbb N_0$-valued random variable with
distribution~$\underline x$. Then, $G_{\mathrm{sel}}
h_{\xi,m}(\underline x) = - \alpha  \operatorname{Cov}[Z, e^{\xi
(Z\wedge m)}] \leq0,$ since two increasing transformations of a
random variable $Z$ have a nonnegative correlation, or, in other
words, the singleton family $\{Z\}$ is associated.)

In the next step we prove \eqref{eqexpBounds0}. We write
%
%e6 ###
\renewcommand{\theequation}{\arabic{section}.\arabic{equation}}
\begin{eqnarray}
\label{eqhxim}
\frac{d}{dt} \mathbf E[h_{\xi,m}(\underline X(t))] & = &\mathbf
E[ G^\alpha_{\mathcal X} h_{\xi,m}(\underline X(t))]
\leq\mathbf E[ G_{\mathrm{mut}} h_{\xi,m}(\underline
X(t))] \nonumber\\
& =& \lambda(e^{\xi} -1) \cdot\mathbf E\Biggl[ \sum_{k=0}^{m-1}
X_k(t) e^{\xi k}\Biggr] \\
& \leq&\lambda(e^\xi- 1) \mathbf
E[h_{\xi,m}(\underline X(t))].\nonumber
\end{eqnarray}
So, by Gronwall's inequality,
\[
\mathbf E[h_{\xi,m}(\underline X(t))] \leq h_{\xi,m}(\underline x)
\cdot\exp\bigl( \lambda t(e^\xi- 1)\bigr)
\]
which gives \eqref{eqexpBounds0} by monotone convergence.

Finally, by Doob's submartingale inequality and monotone
convergence, using~\eqref{eqhxim},
\begin{eqnarray*}
&&\mathbf P\Bigl[  \sup_{0\leq t\leq T} h_\xi(\underline
X(t))>C\Bigr] \\[2pt]
&&\qquad= \lim_{m\to\infty}\mathbf P\Bigl[\sup_{0\leq t\leq T}
h_{\xi,m}(\underline X(t))>C\Bigr] \\[2pt]
 &&\qquad \leq
\lim_{m\to\infty}\mathbf P\biggl[\sup_{0\leq t\leq T} \biggl(
h_{\xi,m}(\underline X(t))\\[2pt]
&&\hspace*{78pt}\qquad{} - \int_0^t
G_{\mathrm{cmut}}h_{\xi,m}(\underline X(s))       +
G_{\mathrm{sel}}h_{\xi,m}(\underline X(s))\,ds\biggr) >C\biggr] \\[2pt]
&&\qquad \leq
\frac1C \lim_{m\to\infty} \mathbf E\biggl[h_{\xi,m}(\underline
X(T))- \int_0^T G_{\mathrm{cmut}}h_{\xi,m}(\underline X(s)) +
G_{\mathrm{sel}}h_{\xi,m}(\underline X(s))\,ds \biggr] \\[2pt]
&&\qquad \leq
\frac1C \lim_{m\to\infty} \biggl(h_{\xi,m}(\underline x) + \int_0^T
\mathbf E[ G_{\mathrm{mut}}h_{\xi,m}(\underline X(s))]
\,ds\biggr) \\[2pt]
&&\qquad \leq\frac1C \biggl( h_\xi(\underline x) +
\lambda(e^\xi-1) h_\xi(\underline x) \int_0^T \exp\bigl(\lambda
s(e^\xi-1)\bigr)\,ds\biggr),
\end{eqnarray*}
and the result follows.
\end{pf}
For the change of measure applied in the next subsection, we will need
that the martingale property of $N_f$ extends from $\mathcal F$ to a
wider class of functions.
\begin{lemma}
Let \label{lext} $\underline x \in\mathbb S_\xi$ for some $\xi>0$
and $\mathcal X = (\underline X_t)_{t\geq0}$ be a solution of the
$(\mathbb S, {\underline x}, G_{\mathcal X}^\alpha, \mathcal
F)$-martingale problem and
%
%e7 ###
\renewcommand{\theequation}{\arabic{section}.\arabic{equation}}
\begin{eqnarray}
\label{eqfvarphi}
f = f_{\varphi_1,\ldots,\varphi_n} \mbox{ be of the form
\eqref{eq122} with }|\varphi_i(\cdot)|\leq C e^{\zeta\cdot}
\nonumber
\\[-8pt]
\\[-8pt]
\eqntext{\mbox{for some }C>0\mbox{ and }\zeta<\xi, i=1,\ldots,n.}
\end{eqnarray}
Then $(N_f(t))_{t\geq0}$, given by~\eqref{eqNf}, is a martingale.
\end{lemma}

\begin{pf} We first observe that $G_{\mathcal X}^\alpha
f(\underline X(t))$ exists for all $t \ge0$; hence~$N_f$
is~well defined. For $\varphi\in\mathcal M(\mathbb N_+)$, let
$\varphi^m(k) := \varphi(k\wedge m)$. We note that\break
$\sum_{k=0}^\infty x_k \varphi^m(k) = \varphi(m) + \sum_{k=0}^m x_k
(\varphi(k) - \varphi(m))$ for $\underline x\in\mathbb S$. Hence,
for $f_{\varphi_1,\ldots,\varphi_n}$ as given in the lemma, the
function $f_{\varphi_1^m,\ldots,\varphi_n^m}$ coincides on $\mathbb S$
with a function in $\mathcal F$. Clearly,
$(N_{f_{\varphi_1^m,\ldots,\varphi_n^m}}(t))_{t\geq0}$ is a martingale
by assumption for all $m=0,1,2,\ldots$ Using~\eqref{eqexpBounds0} and
dominated convergence,
\begin{eqnarray*}
\mathbf E[N_{f_{\varphi_1,\ldots,\varphi_n}}(t)|(\underline
X(r))_{r\leq s}] & =& \lim_{m\to\infty} \mathbf
E[N_{f_{\varphi_1^m,\ldots,\varphi_n^m}}(t)|(\underline X(r))_{r\leq
s}] \\
& =& \lim_{m\to\infty}
N_{f_{\varphi_1^m,\ldots,\varphi_n^m}}(s) =
N_{f_{\varphi_1,\ldots,\varphi_n}}(s).
\end{eqnarray*}
In other words, $(N_f(t))_{t\geq0}$ is a martingale.
\end{pf}

%s4.2 ###
\subsection{Girsanov change of measure}
\label{ssgirs}
In Proposition \ref{P2} we establish a change of measure which shifts
the selection coefficient $\alpha$ of Muller's ratchet with
compensatory mutations. Two assertions from semimartingale theory
which will be required in the proof are recalled in the next remark.

\begin{remark}\label{remmart1}
(1) \textit{A condition for a local martingale to be a martingale}:
Let $\mathcal N = (N_t)_{t\geq0}$ be a local martingale. If
$\mathbf E[\sup_{0\leq t\leq T}|N_t|]<\infty$ for all $T>0$, then~$\mathcal N$
is a martingale; see, for example,\ \citet{Protter2004},
Theorem I.51.\vspace*{-6pt}
\begin{longlist}[(2)]
\item[(2)] \textit{Girsanov Theorem for continuous semimartingales}: Let
$\mathcal L=(L_t)_{t\geq0}$ be a continuous $\mathbf
P$-martingale for some probability measure $\mathbf P$ such that
$\mathcal Z = (Z_t)_{t\geq0}$, given by $Z_t = e^{L_t - 1/2
\langle\mathcal L\rangle_t}$, is a martingale (where $\langle
\mathcal L\rangle$ is the predictable quadratic variation of
$\mathcal L$). If $\mathcal N = (N_t)_{t\geq0}$ is a $\mathbf
P$-local martingale, and $\mathbf Q$ is defined via
\[
\frac{d\mathbf Q}{d\mathbf P}\bigg|_{\mathcal F_t} = Z_t,
\]
then $\mathcal N - \langle\mathcal L, \mathcal N\rangle$ is a
$\mathbf Q$-local martingale. Here, $\langle\mathcal L, \mathcal
N\rangle$ is the (predictable) covariation process between
$\mathcal L$ and $\mathcal N$; see, for example, \citet{Kallenberg2002}, Theorem~18.19
and Lemma~18.21.
\end{longlist}
\end{remark}

\begin{proposition}[(Change of measure)]
For $\underline y\in\mathbb S$, let
\renewcommand{\theequation}{\arabic{section}.\arabic{equation}}
\begin{equation}
\label{eqkappa12}
\kappa_1(\underline y) := \sum_{k=0}^\infty k y_k,\qquad
\kappa_2(\underline y) := \sum_{k=0}^\infty\bigl(k -
\kappa_1(\underline y)\bigr)^2 y_k
\end{equation}
be the expectation and variance of $\underline y$, provided they
exist. Let $\underline x \in\mathbb S_\xi$ for some $\xi>0$ and
$\mathcal X = (\underline X(t))_{t\geq0}$ be a solution of the
$(\mathbb S, {\underline x}, G_{\mathcal X}^\alpha, {\mathcal
F})$-martingale problem, and denote its distribution by $\mathbf
P^\alpha$. Then, the process $\mathcal Z^{\alpha,\alpha'} =
(Z^{\alpha,\alpha'}_t)_{t\geq0}$, given by
%
%e8 ###
\begin{eqnarray}
\label{eqZt}
\qquad Z^{\alpha,\alpha'}_t  &=& \exp\biggl( N(\alpha- \alpha')
\biggl(\kappa_1(\underline X(t)) - \kappa_1(\underline x)
\nonumber
\\[-8pt]
\\[-8pt]
\nonumber
&&\hspace*{73pt}{}
   - \int_0^t \lambda-\gamma
\kappa_1(\underline X(s)) -\frac{\alpha+\alpha'}2
\kappa_2(\underline X(s)) \,ds \biggr)\biggr)
\end{eqnarray}
is a $\mathbf P^\alpha$-local martingale. If $\alpha'>\alpha$, it is
even a $\mathbf P^\alpha$-martingale, and the probability measure
$\mathbf P^{\alpha'}$, defined by
\[
\frac{d \mathbf P^{\alpha'}}{d \mathbf P^\alpha}\bigg|_{\mathcal F_t} =
Z^{\alpha,\alpha'}_t
\]
solves the $(\mathbb S, {\underline x}, G_{\mathcal X}^{\alpha'},
{\mathcal F})$-martingale problem.\vadjust{\goodbreak}
\end{proposition}

\begin{pf}
The proof is an application of the Girsanov transform for continuous
semimartingales; see Remark~\ref{remmart1}.2. By assumption, the
process $\mathcal X$ is continuous, and so is the processes
$(f(\underline X(t)))_{t\geq0}$ for $f$ as in~\eqref{eqfvarphi}.
Set
\[
g(\underline x):=N(\alpha- \alpha')\kappa_1(\underline x),
\]
and define
$\mathcal L = (L_t)_{t\geq0}$ by
\begin{eqnarray*}
L_t & =& N(\alpha- \alpha') \biggl(\kappa_1(\underline X(t)) -
\kappa_1(\underline X(0)) - \int_0^t  G_{\mathcal
X}^\alpha
\kappa_1(\underline X(s)) \,ds\biggr) \\
& =& N(\alpha- \alpha')
\biggl(\kappa_1(\underline X(t)) - \kappa_1(\underline X(0))-
\int_0^t  \lambda- \gamma\kappa_1(\underline
X(s)) - \alpha
\kappa_2(\underline X(s)) \,ds\biggr).
\end{eqnarray*}
Then, $\mathcal L$ is a $\mathbf P^\alpha$-martingale by
Lemma~\ref{lext} with quadratic variation
\begin{eqnarray*}
\langle\mathcal L\rangle_t & =& N^2(\alpha-\alpha')^2 \int_0^t
G^N_{\mathrm{res}} ( \kappa_1(\underline X(s)))^2 \,ds \\
& =&
N(\alpha-\alpha')^2 \int_0^t \kappa_2(\underline X(s))\,ds.
\end{eqnarray*}
For $f\in{\mathcal F}$, let $\mathcal N_f = (N_f(t))_{t\geq0}$ be
as in~\eqref{eqNf}. Then, for $f=f_\varphi\in{\mathcal F}$,
\begin{eqnarray*}
\langle\mathcal L, \mathcal N^f\rangle_t & =& \int_0^t
G^N_{\mathrm{res}} (g(\underline X(s)) f(\underline X(s))) -
g(\underline X(s)) G^N_{\mathrm{res}} f(\underline X(s))\,ds \\
& =&
\frac{\alpha-\alpha'}2 \int_0^t \sum_{k,\ell=0}^\infty
X_k(s)\bigl(\delta_{k\ell} - X_\ell(s)\bigr)\bigl(\varphi(k) \ell+
k\varphi(\ell)\bigr)\,ds \\
 & = &\int_0^t G_{\mathrm{sel}}^{\alpha'}
f(\underline X(s)) - G_{\mathrm{sel}}^{\alpha} f(\underline X(s))
\,ds.
\end{eqnarray*}
By an analogous calculation, one checks that the same identity is
valid for all $f\in{\mathcal F}$. Since $\mathcal L$ is a $\mathbf
P^\alpha$(-local) martingale, the process $\mathcal Z^{\alpha,
\alpha'}$ as well is a~$\mathbf P^\alpha$-local martingale; see~\citet
{Kallenberg2002}, Lemma~18.21.

If $\alpha<\alpha'$ and $\underline x \in\mathbb S_\xi$ [and since
$e^{\xi\kappa_1(\underline x)} \leq h_\xi(\underline x)$ by
Jensen's inequality], we have that $\mathbf E[\sup_{0\leq t\leq T}
Z^{\alpha,\alpha'}_t]<\infty$. Hence, using
Remark~\ref{remmart1}.1, we see that $\mathcal Z^{\alpha,
\alpha'}$\vadjust{\goodbreak}
is a $\mathbf P^\alpha$-martingale. The above calculations and the
Girsanov theorem for continuous semimartingales (recalled in Remark
 \ref{remmart1}.2) then show that
\[
N_f(t) - \langle\mathcal L, \mathcal N_f\rangle_t = f(\underline
X(t)) - \int_0^t G_{\mathcal X}^{\alpha'} f(\underline X(s))\,ds
\]
is a $\mathbf P^{\alpha'}$-martingale. Since $f \in{\mathcal F}$
was arbitrary, $\mathbf P^{\alpha'}$ solves the $(\mathbb S,
{\underline x}, G_{\mathcal X}^{\alpha'}, {\mathcal F})$-martingale
problem.
\end{pf}

%s4.3 ###
\subsection{\texorpdfstring{Proof of Theorem~\protect\ref{T1}}{Proof of Theorem 1}}\label{SproofT1}
First we will prove Proposition~\ref{P34} on the
well-posedness of the martingale problem for $G^\alpha_{\mathcal
X}$. The proof of Theorem~\ref{T1} will then be completed by
observing that a process solves the system of SDEs~\eqref{eqSDEs} if
and only if
it solves the martingale problem for $G_{\mathcal X}^\alpha$ (Lemma
\ref{equiv}).

\begin{pf*}{Proof of Proposition \protect\ref{P34}}
\textit{Step 1: Existence of a solution of the martingale problem}:
For $\alpha=0$, it follows from classical theory
[e.g., \citet{Dawson1993}, Theorem 5.4.1] that the $(\mathbb S,
{\underline x}, G_{\mathcal X}^0, {\mathcal F})$-martingale problem
has a unique solution~$\mathbf P^0$. By Proposition~\ref{P2}, the
change of measure using the martingale $\mathcal Z^{0,\alpha}$ leads
to a distribution $\mathbf P^\alpha$ that solves the $(\mathbb S,
\underline x, G_{\mathcal X}^\alpha, {\mathcal F})$-martingale
problem. This establishes existence.\vspace*{9pt}

\textit{Step 2: Uniqueness of solutions of the martingale problem}:
As in Step 1, let $\mathbf P^0$ be the unique solution of the
$(\mathbb S, {\underline x}, G_{\mathcal X}^0, {\mathcal
F})$-martingale problem. Assume $\mathbf P_1^\alpha$ and $\mathbf
P_2^\alpha$ are two different solutions of the $(\mathbb S,
{\underline x}, G_{\mathcal X}^\alpha, {\mathcal F})$-martingale
problem. Let $\tau_1, \tau_2,\ldots$ be stopping times with
$\tau_n\to\infty$ as $n\to\infty$ such that $(Z^{\alpha,0}_{t\wedge
\tau_n})_{t\geq0}$, given by~\eqref{eqZt}, is both a $\mathbf
P_1^\alpha$-martingale and a $\mathbf P_2^\alpha$-martingale. Since
$\mathbf P_1^\alpha\neq\mathbf P_2^\alpha$, there must be $t\geq0$
such that the distributions of $\underline X(t)$ under $\mathbf
P_1^\alpha$ and $\mathbf P_2^\alpha$ are different; see Theorem 4.4.2
in \citet{EthierKurtz1986}. Hence, there is an $n \in\mathbb N$
such that the distributions of $\underline X(t\wedge\tau_n)$ under
$\mathbf P_1^\alpha$ and $\mathbf P_2^\alpha$ are different. Since
$Z_{t\wedge\tau_n}^{\alpha,0}$ is positive $\mathbf
P_1^\alpha$-a.s. and $\mathbf P_2^\alpha$-a.s., then also the
distributions of $\underline X({t\wedge\tau_n})$ under $Z_{t\wedge
\tau_n}^{\alpha,0} \cdot\mathbf P_1^\alpha$ and $Z_{t\wedge
\tau_n}^{\alpha,0} \cdot\mathbf P_2^\alpha$ are different.
However, by the same arguments as in the proof of
Proposition~\ref{P2}, $Z_{t\wedge\tau_n}^{\alpha,0} \cdot\mathbf
P_1^\alpha$ as well as $Z_{t\wedge\tau_n}^{\alpha,0} \cdot\mathbf
P_2^\alpha$ equal $\mathbf P^0$ on the $\sigma$-algebra
$\sigma((\underline X(s))_{0\leq s\leq t\wedge\tau_n})$, which
contradicts the assumed inequality of $\mathbf P_1^\alpha$ and
$\mathbf P_2^\alpha$. Thus, uniqueness of the $(\mathbb S, \underline
x, G^\alpha_{\mathcal X}, {\mathcal F})$-martingale problem follows.
\end{pf*}

\begin{lemma} [(Equivalence of SDEs and martingale
problem)]\label{equiv}
For $\underline x \in\mathbb S_\xi$, a~process $\mathcal X =
(\underline X(t))_{t\geq0}$ is a weak solution of the system of
SDEs~\eqref{eqSDEs} starting in $\underline x$ if and only if the distribution
of $\mathcal X$ is a solution to the $(\mathbb S, {\underline x},
G_{\mathcal X}^\alpha, \mathcal F)$-martingale problem.
\end{lemma}

\begin{pf}
(1) Assume that $\mathcal X = (\underline X(t))_{t\geq0}$ solves the
system of SDEs~\eqref{eqSDEs}. Then, as a direct consequence of
It\^o's lemma, the distribution of $\mathcal X$ is a~solution to the
$(\mathcal S, \underline x, G_{\mathcal X}^\alpha, \mathcal
F)$-martingale problem.

(2) Conversely, let $\mathcal X = (\underline X(t))_{t\geq0}$ solve
the $(\mathbb S, {\underline x}, G_{\mathcal X}^\alpha, \mathcal
F)$-martingale problem. To see that $\mathcal X$ is a weak solution
of \eqref{eqSDEs}, we may appeal to\vadjust{\goodbreak}
\citet{DaPrato1992}, Theorem~8.2. Specifically, in their notation,
choose $H$ as the Hilbert-space of square-summable elements of
$\mathbb R^{\mathbb N_0}$, $M = (\mathcal N_{f^k})_{k=0,1,2,\ldots}$
with $f^k(\underline x):= x_k$, let $Q$ be the identity on $\mathbb
R^{{\mathbb N_0}\choose 2}= \mathbb\{(w_{k\ell})_{k<\ell}\dvtx
w_{k\ell}\in\mathbb R\}$ and let $\Phi(s)\dvtx  \mathbb
R^{{\mathbb N_0}\choose 2}\to\mathbb R^{\mathbb N_0}$ be given through
the matrix $\Phi(s)_{i,k\ell} := (\delta_{ik} - \delta_{i\ell})
\sqrt{X_k(s) X_\ell(s)}$.
\end{pf}
%

%s5 ###
\section{\texorpdfstring{Proof of Theorem~\protect\ref{T2}}{Proof of Theorem 2}}
\label{Sproof2}
The key element in the proof of Theorem~\ref{T2} is
Proposition~\ref{P2} which represents a solution of~\eqref{eqode1}
through a Markov jump process. For uniqueness of the solution we rely
on a duality derived in Section~\ref{sscalc}. The proof of
Theorem~\ref{T2} is given in Section~\ref{ssproof2}.

%s5.1 ###
\subsection{A particle system}
\label{sspart}
As a preparation to the proof of Theorem~\ref{T2}, we represent the
system of ordinary differential equations by a jump process
$(K_t)_{t\geq0}$. Almost surely, the process will be killed (i.e.,
hit a cemetery state) in finite time. We show in Proposition~\ref{P2}
that a solution of~\eqref{eqode1} is given by the distribution of~$K_t$ conditioned on not being killed by time $t$, $t \ge0$.

\begin{definition}[(Jump process)]
\label{defjump}
Let $(K_t)_{t\geq0}$ be a pure Markov jump process which takes
values in $\{\dagger, 0,1,2,\ldots\}$ and jumps from $k$ to $k+1$ at
rate $\lambda$, from $k$ to $k-1$ at rate $k\gamma$, and from $k$ to
the cemetery state $\dagger$ with rate $\alpha k$.
\end{definition}

\begin{proposition}[(Particle representation)]\label{P2}
Let $\underline x(0) \in\mathbb S_\xi$ for some \mbox{$\xi>0$} and
$(K_t)_{t\geq0}$ be as in Definition~\ref{defjump} with initial
distribution given by $\mathbf P[K_0=k] = x_k(0)$. Then
\renewcommand{\theequation}{\arabic{section}.\arabic{equation}}
\begin{equation}\label{eqoderepr}
x_k(t) := \mathbf P[K_t=k|K_t\neq\dagger]
\end{equation}
solves system \eqref{eqode1}.
\end{proposition}

\begin{pf}
From the definition of $(K_t)_{t\geq0}$, it is clear that for small
$\varepsilon>0$,
\begin{eqnarray*}
&& x_k (t+\varepsilon)\\
&&\qquad= \frac{x_k(t)(1-\alpha k
\varepsilon) + \lambda(x_{k-1}(t) - x_k(t)) \varepsilon+
\gamma((k+1)x_{k+1}(t) - kx_k(t)) \varepsilon}{1-\alpha
\sum_{\ell=0}^\infty\ell x_\ell(t) \varepsilon}  \\
&&\qquad\quad{}+ \mathcal O(\varepsilon^2)\\
&&\qquad = x_k(t) + \Biggl(
-\alpha\Biggl(k - \sum_{\ell=0}^\infty\ell x_\ell(t)\Biggr) x_k(t)
+ \lambda\bigl(x_{k-1}(t) - x_k(t)\bigr) \\
&&\hspace*{148pt}{}
  + \gamma\bigl((k+1)x_{k+1}(t) - k x_k(t)\bigr)\Biggr)\varepsilon
+ \mathcal O(\varepsilon^2),
\end{eqnarray*}
which implies the result as $\varepsilon\to0$.
\end{pf}

%s5.2 ###
\subsection{Dynamics of the cumulant generating function}
\label{sscalc}
The proof of uniqueness of~\eqref{eqode1} requires some preliminary
computations which we carry out next. Recall the function $h_\xi$
from~\eqref{eqhdef}. Note the function $\zeta\mapsto\log
h_\zeta(\underline x)$ is the cumulant generating function of
$\underline x \in\mathbb S$. Cumulants have already been proven to be
useful in the study of Muller's ratchet; see
\citet{EtheridgePfaffelhuberWakolbinger2007}. Here, we compute the
dynamics of the cumulant generating function.

\begin{proposition}[(Dynamics of cumulant generating
function)]\label{P43}
For any solution $t\mapsto\underline x(t)$ of \eqref{eqode1}
taking values in $\mathbb S_\xi$ for $\xi>0$ and $0<\zeta<\xi$,
\[
\frac{d}{dt} \log h_\zeta(\underline x(t)) = \alpha\sum_{\ell
=0}^\infty
\ell x_\ell(t) + \lambda(e^\zeta-1) - \bigl(\alpha+
\gamma(1-e^{-\zeta})\bigr)\frac{d}{d\zeta} \log h_\zeta(\underline x(t)).
\]
\end{proposition}

\begin{pf}
Abbreviating $\underline x := \underline x(t)$, we compute
\begin{eqnarray*}
h_\zeta(\underline x) \frac{d}{dt}  \log h_\zeta(\underline x)& =&
\alpha\sum_{\ell=0}^\infty\sum_{k=0}^\infty(\ell-k)x_\ell x_k
e^{\zeta k} + \lambda\sum_{k=0}^\infty(x_{k-1} -x_k)e^{\zeta k}
\\
&  &{}      + \gamma
\sum_{k=0}^\infty\bigl((k+1)x_{k+1} - kx_k\bigr)e^{\zeta k} \\
& = &\alpha
\Biggl( \Biggl(\sum_{\ell=0}^\infty\ell x_\ell\Biggr)
h_\zeta(\underline x) - \frac{d}{d\zeta} h_\zeta(\underline
x)\Biggr) + \lambda(e^\zeta-1) h_\zeta(\underline x) \\
&&{}  - \gamma
(1-e^{-\zeta}) \frac{d}{d\zeta} h_\zeta(\underline x)
\end{eqnarray*}
and so
\[
\qquad\frac{d}{dt} \log h_\zeta(\underline x)  = \alpha\sum_{\ell
=0}^\infty\ell x_\ell+ \lambda(e^\zeta-1)-\bigl (\alpha+
\gamma(1-e^{-\zeta})\bigr)\frac{d}{d\zeta} \log h_\zeta(\underline
x).\qquad\qed
\]
\noqed\end{pf}

The equation in Proposition \ref{P43} relates the time-derivative of
$\log h_\zeta(\underline x(t))$ with the $\zeta$-derivative of the
same function and leads to a \textit{duality relation} formulated in
Corollary \ref{cordual}. In Markov process theory, dualities are
particularly useful to obtain uniqueness results; cf. \citet
{EthierKurtz1986}, page 188ff.
Our application in Section
\ref{ssproof2} will be in this spirit.
\begin{corollary}[(Duality)]\label{cordual}
Let $t\mapsto\underline x(t)$ be a
solution of~\eqref{eqode1} taking values in $\mathbb S_\xi$ for
$\xi>0$. Moreover let $\zeta\dvtx  t\mapsto\zeta(t)$ be the solution of
$\zeta' = -(\alpha+ \gamma(1-e^{-\zeta}))$, starting in some
$\zeta(0)<\xi$. Then
\[
\log h_{\zeta(0)}(\underline x(t)) = \log h_{\zeta(t)}(\underline
x(0)) + \int_0^t \Biggl(\lambda\bigl(e^{\zeta(t-s)}-1\bigr) + \sum_{\ell
=0}^\infty\ell x_\ell(s)\Biggr) \,ds.
\]
\end{corollary}

\begin{pf}
Using Proposition~\ref{P43} and noting, for any differentiable $g:
\zeta\mapsto g(\zeta)$, the equality
\[
\frac{d}{ds} g\bigl(\zeta(t-s)\bigr) = \bigl(\alpha+ \gamma\bigl(1-e^{-\zeta(t-s)}\bigr)\bigr)\frac
{d}{d\zeta}
g\bigl(\zeta(t-s)\bigr),
\]
we obtain
\[
\frac{d}{ds} \log h_{\zeta(t-s)}(\underline x(s)) = \lambda
\bigl(e^{\zeta(t-s)}-1\bigr) + \alpha\sum_{\ell=0}^\infty\ell x_\ell(s).
\]
Now the assertion follows by integrating.
\end{pf}

%s5.3 ###
\subsection{\texorpdfstring{Proof of Theorem~\protect\ref{T2}}{Proof of Theorem 2}}
\label{ssproof2}
We proceed in two steps. First, we derive an explicit solution
of~\eqref{eqode1} by using Proposition~\ref{P2}, that is, by computing
the distribution of the jump process $(K_t)_{t\geq0}$ conditioned on
not being killed by time $t$. This will result in the right-hand side
of~\eqref{eqT11}. In a second step, we show uniqueness of solutions
of~\eqref{eqode1} in $\mathbb S_\xi$.

\textit{Step 1: Computation of the right-hand side of}
\eqref{eqoderepr}:
In order to derive an explicit formula for the probability specified
in \eqref{eqoderepr}, we note that the process $(K_t)_{t\geq0}$ can
be realized as the following mutation-couting process:
\begin{itemize}
\item Start with $K_0$ mutations, with the random number $K_0$
distributed according to $(x_k(0))_{k=0,1,2,\ldots}$.
\item New mutations arise at rate $\lambda$.
\item Every mutation (present from the start or newly arisen) starts
an exponential waiting time with parameter $\alpha+\gamma$. If this
waiting time expires, then with probability
$\frac{\alpha}{\alpha+\gamma}$ the process jumps to $\dagger$, and
with the complementary probability $\frac{\gamma}{\alpha+\gamma}$
the mutation disappears.
\end{itemize}
With $x_k(t)$ defined by \eqref{eqoderepr}, we decompose the
probability of the event $\{K_t = k\}$ with respect to the number of
mutations present at time $0$. If $K_0=i$, a number $j\leq i \wedge k$
of these initial mutations are not compensated by time $t$, and the
remaining $i-j$ are compensated. In addition, a number $l\geq k-j$
mutations arise at times $0\leq t_1\leq\cdots\leq t_{l}\leq t$. From
these, \mbox{$l-k+j$} are compensated, and the remaining $k-j$ are not
compensated. These arguments lead to the following calculation, where
we write $\sim$ for equality up to factors not depending on $k$. The
first $\sim$ comes from the fact that the right-hand side is the
unconditional probability $\mathbf P[K_t=k]$,
\begin{eqnarray*}
x_k(t) & \sim&\sum_{i=0}^\infty x_i(0) \sum_{j=0}^{i\wedge k}
\pmatrix{{i}\vspace*{2pt}\cr{j}} \biggl(\frac{\gamma}{\alpha+\gamma}\bigl(
1-e^{-(\alpha+\gamma)t}\bigr)\biggr)^{i-j} \cdot e^{-j(\alpha+\gamma)t}
\\[2pt]
&  &{}\times \!\sum_{l=k-j}^\infty \mathop{ \!\int _{\{\mathcal T =
(t_1,\ldots,t_{l})\dvtx }}_{0\leq t_1\leq\cdots\leq t_{l} \}}\!
d(t_1,\ldots,t_{l}) \lambda^{l} e^{-\lambda t_1}e^{-\lambda(t_2-t_1)}
\cdots e^{-\lambda(t_{l}-t_{l-1})} e^{-\lambda(t-t_{l})} \\[2pt]
&&{}\times
\mathop{\sum_{\mathcal S\subseteq\mathcal T}}_{|\mathcal S| = l-k+j}
\prod_{r\in\mathcal S}\frac{\gamma}{\alpha+\gamma} \bigl(
1-e^{-(\alpha+\gamma)(t-r)}\bigr) \cdot\prod_{s\in\mathcal
T\setminus\mathcal S}e^{-(\alpha+\gamma)(t-s)}
\\[2pt]
& =& \sum_{i=0}^\infty x_i(0) \sum_{j=0}^{i\wedge k}
\pmatrix{{i}\vspace*{2pt}\cr{j}}
\biggl(\frac{\gamma}{\alpha+\gamma}\bigl(
1-e^{-(\alpha+\gamma)t}\bigr)\biggr)^{i-j} \cdot e^{-j(\alpha+\gamma)t}
\\[2pt]
&&{}\times   \sum_{l=k-j}^\infty\frac{\lambda^l}{l!} e^{-\lambda t}
        \mathop{ \int _{\{\mathcal T =
(t_1,\ldots,t_{l})\dvtx }}_{0\leq t_1,\ldots, t_{l} \}}
d(t_1,\ldots,t_{l})\\[2pt]
&&{}\times
\mathop{\sum_{\mathcal S\subseteq\mathcal T}}_{|\mathcal S| = l-k+j}
   \prod_{r\in\mathcal
S}\frac{\gamma}{\alpha+\gamma}\bigl(
1-e^{-(\alpha+\gamma)(t-r)}\bigr) \prod_{s\in\mathcal T\setminus
\mathcal S}e^{-(\alpha+\gamma)(t-s)} \\[2pt]
& \sim&
\sum_{i=0}^\infty x_i(0) \sum_{j=0}^{i\wedge k}
\pmatrix{{i}\vspace*{2pt}\cr{j}}
\biggl(\frac{\gamma}{\alpha+\gamma}\bigl(
1-e^{-(\alpha+\gamma)t}\bigr)\biggr)^{i-j} e^{-j(\alpha+\gamma)t}\\[2pt]
 &&\hspace*{56pt}{}
 \times\frac{\lambda^{k-j}}{(k-j)!}
\biggl(\int_0^te^{-(\alpha+\gamma)(t-s)}\,ds \biggr)^{k-j} \\[2pt]
&&{}
 \times\sum_{l=0}^\infty\frac{\lambda^l}{l!} \biggl(
\frac{\gamma}{\alpha+\gamma}\int_0^t
1-e^{-(\alpha+\gamma)(t-r)}dr\biggr)^{l}
 \\[2pt] & \sim&
\sum_{i=0}^\infty x_i(0) \sum_{j=0}^{i\wedge k}
\pmatrix{{i}\vspace*{2pt}\cr{j}}
\biggl(\frac{\gamma}{\alpha+\gamma}\bigl(
1-e^{-(\alpha+\gamma)t}\bigr)\biggr)^{i-j} e^{-j(\alpha+\gamma)t}\\[2pt]
&&\hspace*{56pt}{}
  \times\frac{\lambda^{k-j}}{(k-j)!}
\biggl(\frac{1}{\alpha+\gamma} \bigl(1-e^{-(\alpha+\gamma)t}\bigr) \biggr)^{k-j},
\end{eqnarray*}
where the first ``$=$'' comes from the symmetry of the
integrand. Summing the right-hand side gives
\begin{eqnarray*}
&&\sum_{i=0}^\infty x_i(0) \sum_{j=0}^i
\pmatrix{{i}\vspace*{2pt}\cr{j}}
\biggl(\frac{\gamma}{\alpha+\gamma}\bigl(
1-e^{-(\alpha+\gamma)t}\bigr)\biggr)^{i-j} e^{-j(\alpha+\gamma)t} \\[2pt]
&&\quad{}   \times\sum_{k=j}^\infty
\frac{\lambda^{k-j}}{(k-j)!} \biggl( \frac{1}{\alpha+\gamma}
\bigl(1-e^{-(\alpha+\gamma)t}\bigr) \biggr)^{k-j} \\[2pt]
 &&\qquad = \sum_{i=0}^\infty
x_i(0) \biggl(\frac{\gamma}{\alpha+\gamma} -
\frac{\alpha}{\alpha+\gamma} e^{-(\alpha+\gamma)t}\biggr)^i \cdot
\exp\biggl(\frac{\lambda}{\alpha+\gamma}\bigl(1-e^{-(\alpha+\gamma)t}\bigr)
\biggr).
\end{eqnarray*}
Hence,
\begin{eqnarray*}
x_k(t) & =& \Biggl( \sum _{i=0}^\infty x_i(0)
\sum _{j=0}^{i\wedge k} \pmatrix{{i}\cr{j}} \bigl(\bigl({\gamma\bigl(
1-e^{-(\alpha+\gamma)t}\bigr)}\bigr)/({\alpha+\gamma})\bigr)^{i-j}\\
&&\hspace*{40pt}\qquad{}\times e^{-j(\alpha+\gamma)t} \bigl({1}/{(k-j)!}\bigr)
\\
&&\hspace*{53pt}\qquad{}\times \bigl({\bigl(\lambda\bigl(1-e^{-(\alpha+\gamma)t}\bigr)\bigr)}/{(\alpha+\gamma)}
\bigr)^{k-j}\Biggr)\\
&&\Bigg/\Biggl( \sum _{i=0}^\infty x_i(0)
\bigl({\gamma}/{(\alpha+\gamma)} - {\alpha}/{(\alpha+\gamma)}
e^{-(\alpha+\gamma)t}\bigr)^i \\
&&\hspace*{46pt}\qquad{}\times
\exp\bigl({\lambda}/{(\alpha+\gamma)}\bigl(1-e^{-(\alpha+\gamma)t}\bigr)
\bigr)\Biggr)
\end{eqnarray*}
which shows \eqref{eqT11}. To see \eqref{eqT12}, it suffices to note
that all terms in the numerator, except for $j=0$, converge to $0$
as $t\to\infty$. Hence the result is proved.

\textit{Step 2: Uniqueness in $\mathbb S_\xi$}:
Let $(\underline y(t))_{t\geq0}$ be a solution of~\eqref{eqode1}
starting in $\underline y(0)\in\mathbb S_\xi$. From the analog of
Lemma~\ref{lexpBounds} in the case $N=\infty$ we have
$h_\xi(\underline y(t)) \leq h_\xi(\underline y(0))\exp(\lambda
t(e^\xi-1)) < \infty$, that is, $\underline y(t) \in\mathbb S_\xi$ for
all $t\geq0$.

If $(\underline x(t))_{t\geq0}$ and $(\underline y(t))_{t\geq0}$ are
solutions of ~\eqref{eqode1} with $\underline x(0) = \underline y(0)
\in\mathbb S_\xi$, then we obtain from Corollary~\ref{cordual} that
for all $0<\zeta<\xi$ and any $t\geq0$,
\renewcommand{\theequation}{\arabic{section}.\arabic{equation}}
\begin{equation}
\label{equnique1}
\log h_{\zeta}(\underline x(t)) - \log h_{\zeta}(\underline y(t)) =
\int_0^t \sum_{\ell=0}^\infty\ell\bigl(x_\ell(s) - y_\ell(s)\bigr) \,ds.
\end{equation}
Since only the left-hand side depends on $\zeta$, this enforces that
both sides vanish. Indeed, taking derivatives with respect to $\zeta$
at $\zeta=0$ the previous equality gives
\[
\sum_{\ell=0}^\infty\ell\bigl(x_\ell(t) - y_\ell(t)\bigr) = 0,  \qquad t\ge0.
\]
Plugging this back into~\eqref{equnique1} gives
\begin{equation}
\label{equnique2}
\log h_{\zeta}(\underline x(t)) = \log h_{\zeta}(\underline y(t)),
\qquad  t\ge0.
\end{equation}
Since the function $\zeta\mapsto\log h_\zeta(\underline x)$ (for
$0<\zeta<\xi$) characterizes $\underline x \in\mathbb S_\xi$, we
obtain that $\underline x(t) = \underline y(t)$. This completes the
proof of Theorem~\ref{T2}.

\section*{Acknowledgments}
We thank Etienne Pardoux for a helpful discussion, and a referee for
various suggestions that led to an improvement of the presentation. Part
of this work was carried out while P. Pfaffelhuber and A.~Wakolbinger were visiting the program
in Probability and Discrete Mathematics in Mathematical Biology at the
University of Singapore in April 2011, whose hospitality is gratefully
acknowledged.

% imsref loaded by akundreckaite, 2012-04-26 13:19:33

%

%suskaldyti doi

\printaddresses


\begin{thebibliography}{39}
% BibTex style file: ims.bst, 2011-05-30
% Default style options (sort=0,type=number).
% Used options (sort=1,type=nameyear).

%b1 ###
\bibitem[\protect\citeauthoryear{Antezana and
  Hudson}{1997}]{AntezanaHudson1997}
\begin{barticle}[auto:STB|2012/04/25|11:03:30]
\bauthor{\bsnm{Antezana},~\bfnm{M.~A.}\binits{M.~A.}} \AND
  \bauthor{\bsnm{Hudson},~\bfnm{R.~R.}\binits{R.~R.}}
(\byear{1997}).
\btitle{Era reversibile! Point-mutations, the Ratchet, and the initial success
  of eukaryotic sex: A simulation study}.
\bjournal{Evolutionary Theory}
\bvolume{11}
\bpages{209--235}.
\bptok{imsref}%
\end{barticle}
\endbibitem

%b2 ###
\bibitem[\protect\citeauthoryear{Audiffren}{2011}]{Audiffren2011}
\begin{bmisc}[auto:STB|2012/04/25|11:03:30]
\bauthor{\bsnm{Audiffren},~\bfnm{J.}\binits{J.}}
(\byear{2011}).
\bhowpublished{Ph.D. thesis. Univ. de
  Provence, Marseille.}
\bptok{imsref}%
\end{bmisc}
\endbibitem

%b3 ###
\bibitem[\protect\citeauthoryear{Audiffren and
  Pardoux}{2011}]{AudiffrenPardoux2011}
\begin{bmisc}[auto:STB|2012/04/25|11:03:30]
\bauthor{\bsnm{Audiffren},~\bfnm{J.}\binits{J.}} \AND
  \bauthor{\bsnm{Pardoux},~\bfnm{E.}\binits{E.}}
(\byear{2011}).
\bhowpublished{Muller's ratchet clicks in finite time.
Unpublished manuscript.}
\bptok{imsref}%
\end{bmisc}
\endbibitem

%b4 ###
\bibitem[\protect\citeauthoryear{Chao}{1990}]{pmid2247152}
\begin{barticle}[pbm]
\bauthor{\bsnm{Chao},~\bfnm{L.}\binits{L.}}
(\byear{1990}).
\btitle{Fitness of RNA virus decreased by Muller's ratchet}.
\bjournal{Nature}
\bvolume{348}
\bpages{454--455}.
\bid{doi={10.1038/348454a0}, issn={0028-0836}, pmid={2247152}}
\bptok{imsref}%
\end{barticle}
\endbibitem

%b5 ###
\bibitem[\protect\citeauthoryear{Cuthbertson}{2007}]{Cuthbertson2007}
\begin{bmisc}[auto:STB|2012/04/25|11:03:30]
\bauthor{\bsnm{Cuthbertson},~\bfnm{C.}\binits{C.}}
(\byear{2007}).
\bhowpublished{Limits to the rate of adaptation. Ph.D. thesis, Univ. Oxford}.
\bptok{imsref}%
\end{bmisc}
\endbibitem

%b6 ###
\bibitem[\protect\citeauthoryear{Da~Prato and Zabczyk}{1992}]{DaPrato1992}
\begin{bbook}[mr]
\bauthor{\bsnm{Da~Prato},~\bfnm{Giuseppe}\binits{G.}} \AND
  \bauthor{\bsnm{Zabczyk},~\bfnm{Jerzy}\binits{J.}}
(\byear{1992}).
\btitle{Stochastic Equations in Infinite Dimensions}.
\bseries{Encyclopedia of Mathematics and Its Applications}
\bvolume{44}.
\bpublisher{Cambridge Univ. Press}, \baddress{Cambridge}.
\bid{doi={10.1017/CBO9780511666223}, mr={1207136}}
\bptok{imsref}%
\end{bbook}
\endbibitem

%b7 ###
\bibitem[\protect\citeauthoryear{Dawson}{1993}]{Dawson1993}
\begin{bincollection}[mr]
\bauthor{\bsnm{Dawson},~\bfnm{Donald~A.}\binits{D.~A.}}
(\byear{1993}).
\btitle{Measure-valued {M}arkov processes}.
In \bbooktitle{\'{E}cole D'\'{E}t\'e de {P}robabilit\'es de {S}aint-{F}lour
  {XXI}---1991}
(\beditor{\binits{P.} \bsnm{Hennequin}}, ed.).
\bseries{Lecture Notes in Math.}
\bvolume{1541}
\bpages{1--260}.
\bpublisher{Springer}, \baddress{Berlin}.
\bid{doi={10.1007/BFb0084190}, mr={1242575}}
\bptok{imsref}%
\end{bincollection}
\endbibitem

%b8 ###
\bibitem[\protect\citeauthoryear{Dawson and Li}{2010}]{DawsonLi2010}
\begin{barticle}[auto:STB|2012/04/25|11:03:30]
\bauthor{\bsnm{Dawson},~\bfnm{D.}\binits{D.}} \AND
  \bauthor{\bsnm{Li},~\bfnm{Z.}\binits{Z.}}
(\byear{2010}).
\btitle{Stochastic equations, flows and measure-valued processes}.
\bjournal{Ann. Probab.}
\bvolume{40}
\bpages{813--857}.
\bptok{imsref}%
\end{barticle}
\endbibitem

%b9 ###
\bibitem[\protect\citeauthoryear{Desai and Fisher}{2007}]{DesaiFisher2007}
\begin{barticle}[pbm]
\bauthor{\bsnm{Desai},~\bfnm{Michael~M.}\binits{M.~M.}} \AND
  \bauthor{\bsnm{Fisher},~\bfnm{Daniel~S.}\binits{D.~S.}}
(\byear{2007}).
\btitle{Beneficial mutation selection balance and the effect of linkage on
  positive selection}.
\bjournal{Genetics}
\bvolume{176}
\bpages{1759--1798}.
\bid{doi={10.1534/genetics.106.067678}, issn={0016-6731},
  pii={genetics.106.067678}, pmcid={1931526}, pmid={17483432}}
\bptok{imsref}%
\end{barticle}
\endbibitem

%b10 ###
\bibitem[\protect\citeauthoryear{Etheridge, Pfaffelhuber and Wakolbinger}{2009}]{EtheridgePfaffelhuberWakolbinger2007}
\begin{bincollection}[mr]
\bauthor{\bsnm{Etheridge},~\bfnm{Alison~M.}\binits{A.~M.}},
  \bauthor{\bsnm{Pfaffelhuber},~\bfnm{Peter}\binits{P.}} \AND
  \bauthor{\bsnm{Wakolbinger},~\bfnm{Anton}\binits{A.}}
(\byear{2009}).
\btitle{How often does the ratchet click? {F}acts, heuristics, asymptotics}.
In \bbooktitle{Trends in Stochastic Analysis}.
\bseries{London Mathematical Society Lecture Note Series}
\bvolume{353}
\bpages{365--390}.
\bpublisher{Cambridge Univ. Press}, \baddress{Cambridge}.
\bid{mr={2562161}}
\bptnote{check year}%
\bptok{imsref}%
\end{bincollection}
\endbibitem

%b11 ###
\bibitem[\protect\citeauthoryear{Ethier and Kurtz}{1986}]{EthierKurtz1986}
\begin{bbook}[mr]
\bauthor{\bsnm{Ethier},~\bfnm{Stewart~N.}\binits{S.~N.}} \AND
  \bauthor{\bsnm{Kurtz},~\bfnm{Thomas~G.}\binits{T.~G.}}
(\byear{1986}).
\btitle{Markov Processes: Characterization and Convergence}.
%  and Mathematical Statistics}.
\bpublisher{Wiley}, \baddress{New York}.
\bid{doi={10.1002/9780470316658}, mr={0838085}}
\bptok{imsref}%
\end{bbook}
\endbibitem

%b12 ###
\bibitem[\protect\citeauthoryear{Ethier and Kurtz}{1993}]{EthierKurtz1993}
\begin{barticle}[mr]
\bauthor{\bsnm{Ethier},~\bfnm{S.~N.}\binits{S.~N.}} \AND
  \bauthor{\bsnm{Kurtz},~\bfnm{Thomas~G.}\binits{T.~G.}}
(\byear{1993}).
\btitle{Fleming--{V}iot processes in population genetics}.
\bjournal{SIAM J. Control Optim.}
\bvolume{31}
\bpages{345--386}.
\bid{doi={10.1137/0331019}, issn={0363-0129}, mr={1205982}}
\bptok{imsref}%
\end{barticle}
\endbibitem

%b13 ###
\bibitem[\protect\citeauthoryear{Ethier and Shiga}{2000}]{EthierShiga2000}
\begin{barticle}[mr]
\bauthor{\bsnm{Ethier},~\bfnm{Stewart~N.}\binits{S.~N.}} \AND
  \bauthor{\bsnm{Shiga},~\bfnm{Tokuzo}\binits{T.}}
(\byear{2000}).
\btitle{A {F}leming--{V}iot process with unbounded selection}.
\bjournal{J. Math. Kyoto Univ.}
\bvolume{40}
\bpages{337--361}.
\bid{issn={0023-608X}, mr={1787875}}
\bptok{imsref}%
\end{barticle}
\endbibitem

%b14 ###
\bibitem[\protect\citeauthoryear{Gabriel, Lynch and
  B{\"u}rger}{1993}]{GabrielLynchBuerger1993}
\begin{barticle}[auto:STB|2012/04/25|11:03:30]
\bauthor{\bsnm{Gabriel},~\bfnm{W.}\binits{W.}},
  \bauthor{\bsnm{Lynch},~\bfnm{M.}\binits{M.}} \AND
  \bauthor{\bsnm{B{\"u}rger},~\bfnm{R.}\binits{R.}}
(\byear{1993}).
\btitle{Muller's ratchet and mutational meltdowns}.
\bjournal{Evolution}
\bvolume{47}
\bpages{1744--1757}.
\bptok{imsref}%
\end{barticle}
\endbibitem

%b15 ###
\bibitem[\protect\citeauthoryear{Gerrish and Lenski}{1998}]{GerrishLenski98}
\begin{barticle}[auto:STB|2012/04/25|11:03:30]
\bauthor{\bsnm{Gerrish},~\bfnm{P.}\binits{P.}} \AND
  \bauthor{\bsnm{Lenski},~\bfnm{R.}\binits{R.}}
(\byear{1998}).
\btitle{The fate of competing beneficial mutations in an asexual population}.
\bjournal{Genetica}
\bvolume{102/103}
\bpages{127--144}.
\bptok{imsref}%
\end{barticle}
\endbibitem

%b16 ###
\bibitem[\protect\citeauthoryear{Gessler}{1995}]{Gessler1995}
\begin{barticle}[pbm]
\bauthor{\bsnm{Gessler},~\bfnm{D.~D.}\binits{D.~D.}}
(\byear{1995}).
\btitle{The constraints of finite size in asexual populations and the rate of
  the ratchet}.
\bjournal{Genet. Res.}
\bvolume{66}
\bpages{241--253}.
\bid{pmid={16553995}}
\bptok{imsref}%
\end{barticle}
\endbibitem

%b17 ###
\bibitem[\protect\citeauthoryear{Gordo and
  Charlesworth}{2000}]{GordoCharlesworth2000}
\begin{barticle}[pbm]
\bauthor{\bsnm{Gordo},~\bfnm{I.}\binits{I.}} \AND
  \bauthor{\bsnm{Charlesworth},~\bfnm{B.}\binits{B.}}
(\byear{2000}).
\btitle{On the speed of Muller's ratchet}.
\bjournal{Genetics}
\bvolume{156}
\bpages{2137--2140}.
\bid{issn={0016-6731}, pmcid={1461355}, pmid={11187462}}
\bptok{imsref}%
\end{barticle}
\endbibitem

%b18 ###
\bibitem[\protect\citeauthoryear{Haigh}{1978}]{Haigh1978}
\begin{barticle}[mr]
\bauthor{\bsnm{Haigh},~\bfnm{John}\binits{J.}}
(\byear{1978}).
\btitle{The accumulation of deleterious genes in a population---{M}uller's
  ratchet}.
\bjournal{Theoret. Population Biol.}
\bvolume{14}
\bpages{251--267}.
\bid{doi={10.1016/0040-5809(78)90027-8}, issn={0040-5809}, mr={0514423}}
\bptok{imsref}%
\end{barticle}
\endbibitem

%b19 ###
\bibitem[\protect\citeauthoryear{Handel, Regoes and Antia}{2006}]{Handel2006}
\begin{barticle}[pbm]
\bauthor{\bsnm{Handel},~\bfnm{Andreas}\binits{A.}},
  \bauthor{\bsnm{Regoes},~\bfnm{Roland~R.}\binits{R.~R.}} \AND
  \bauthor{\bsnm{Antia},~\bfnm{Rustom}\binits{R.}}
(\byear{2006}).
\btitle{The role of compensatory mutations in the emergence of drug
  resistance}.
\bjournal{PLoS Comput. Biol.}
\bvolume{2}
\bpages{e137}.
\bid{doi={10.1371/journal.pcbi.0020137}, issn={1553-7358},
  pii={05-PLCB-RA-0334R3}, pmcid={1599768}, pmid={17040124}}
\bptok{imsref}%
\end{barticle}
\endbibitem

%b20 ###
\bibitem[\protect\citeauthoryear{Higgs and Woodcock}{1995}]{HiggsWoodcock1995}
\begin{barticle}[auto:STB|2012/04/25|11:03:30]
\bauthor{\bsnm{Higgs},~\bfnm{P.~G.}\binits{P.~G.}} \AND
  \bauthor{\bsnm{Woodcock},~\bfnm{G.}\binits{G.}}
(\byear{1995}).
\btitle{The accumulation of mutations in asexual populations and the structure
  of genealogical trees in the presence of selection}.
\bjournal{J. Math. Biol.}
\bvolume{33}
\bpages{677--702}.
\bptok{imsref}%
\end{barticle}
\endbibitem

%b21 ###
\bibitem[\protect\citeauthoryear{Howe and Denver}{2008}]{Howe2008}
\begin{barticle}[auto:STB|2012/04/25|11:03:30]
\bauthor{\bsnm{Howe},~\bfnm{D.~K.}\binits{D.~K.}} \AND
  \bauthor{\bsnm{Denver},~\bfnm{D.~R.}\binits{D.~R.}}
(\byear{2008}).
\btitle{Muller's ratchet and compensatory mutation in \textit{Caenorhabditis
  briggsae} mitochondrial genome evolution}.
\bjournal{BMC Evol. Biol.}
\bvolume{8}
\bpages{62--62}.
\bptok{imsref}%
\end{barticle}
\endbibitem

%b22 ###
\bibitem[\protect\citeauthoryear{Jain}{2008}]{pmid18689884}
\begin{barticle}[pbm]
\bauthor{\bsnm{Jain},~\bfnm{Kavita}\binits{K.}}
(\byear{2008}).
\btitle{Loss of least-loaded class in asexual populations due to drift and
  epistasis}.
\bjournal{Genetics}
\bvolume{179}
\bpages{2125--2134}.
\bid{doi={10.1534/genetics.108.089136}, issn={0016-6731},
  pii={genetics.108.089136}, pmcid={2516084}, pmid={18689884}}
\bptok{imsref}%
\end{barticle}
\endbibitem

%b23 ###
\bibitem[\protect\citeauthoryear{Kallenberg}{2002}]{Kallenberg2002}
\begin{bbook}[mr]
\bauthor{\bsnm{Kallenberg},~\bfnm{Olav}\binits{O.}}
(\byear{2002}).
\btitle{Foundations of Modern Probability},
\bedition{2nd} ed.
\bpublisher{Springer}, \baddress{New York}.
\bid{mr={1876169}}
\bptok{imsref}%
\end{bbook}
\endbibitem

%b24 ###
\bibitem[\protect\citeauthoryear{Loewe}{2006}]{Loewe2006}
\begin{barticle}[pbm]
\bauthor{\bsnm{Loewe},~\bfnm{Laurence}\binits{L.}}
(\byear{2006}).
\btitle{Quantifying the genomic decay paradox due to Muller's ratchet in human
  mitochondrial DNA}.
\bjournal{Genet. Res.}
\bvolume{87}
\bpages{133--159}.
\bid{pii={S0016672306008123}, doi={10.1017/S0016672306008123}, pmid={16709275}}
\bptok{imsref}%
\end{barticle}
\endbibitem

%b25 ###
\bibitem[\protect\citeauthoryear{Maia, Botelho and
  Fontanari}{2003}]{MaiaBotelhoFontanari2003}
\begin{barticle}[mr]
\bauthor{\bsnm{Maia},~\bfnm{Leonardo~P.}\binits{L.~P.}},
  \bauthor{\bsnm{Botelho},~\bfnm{Daniela~F.}\binits{D.~F.}} \AND
  \bauthor{\bsnm{Fontanari},~\bfnm{Jos{\'e}~F.}\binits{J.~F.}}
(\byear{2003}).
\btitle{Analytical solution of the evolution dynamics on a
  multiplicative-fitness landscape}.
\bjournal{J. Math. Biol.}
\bvolume{47}
\bpages{453--456}.
\bid{doi={10.1007/s00285-003-0208-8}, issn={0303-6812}, mr={2029007}}
\bptok{imsref}%
\end{barticle}
\endbibitem

%b26 ###
\bibitem[\protect\citeauthoryear{Maier
  et~al.}{2008}]{Maier2008BMC-Biol18755031}
\begin{barticle}[auto:STB|2012/04/25|11:03:30]
\bauthor{\bsnm{Maier},~\bfnm{U.}\binits{U.}},
  \bauthor{\bsnm{Bozarth},~\bfnm{A.}\binits{A.}},
  \bauthor{\bsnm{Funk},~\bfnm{H.}\binits{H.}},
  \bauthor{\bsnm{Zauner},~\bfnm{S.}\binits{S.}},
  \bauthor{\bsnm{Rensing},~\bfnm{S.}\binits{S.}},
  \bauthor{\bsnm{Schmitz-Linneweber},~\bfnm{C.}\binits{C.}},
  \bauthor{\bsnm{B{\"o}rner},~\bfnm{T.}\binits{T.}} \AND
  \bauthor{\bsnm{Tillich},~\bfnm{M.}\binits{M.}}
(\byear{2008}).
\btitle{Complex chloroplast rna metabolism: Just debugging the genetic
  programme?}
\bjournal{BMC Biol.}
\bvolume{6}
\bpages{36--36}.
\bptok{imsref}%
\end{barticle}
\endbibitem

%b27 ###
\bibitem[\protect\citeauthoryear{Maisnier-Patin and
  Andersson}{2004}]{Maisnier2004}
\begin{barticle}[pbm]
\bauthor{\bsnm{Maisnier-Patin},~\bfnm{Sophie}\binits{S.}} \AND
  \bauthor{\bsnm{Andersson},~\bfnm{Dan~I.}\binits{D.~I.}}
(\byear{2004}).
\btitle{Adaptation to the deleterious effects of antimicrobial drug resistance
  mutations by compensatory evolution}.
\bjournal{Res. Microbiol.}
\bvolume{155}
\bpages{360--369}.
\bid{doi={10.1016/j.resmic.2004.01.019}, issn={0923-2508},
  pii={S0923-2508(04)00068-3}, pmid={15207868}}
\bptok{imsref}%
\end{barticle}
\endbibitem

%b28 ###
\bibitem[\protect\citeauthoryear{Maynard~Smith}{1978}]{MaynardSmith1978}
\begin{bbook}[auto:STB|2012/04/25|11:03:30]
\bauthor{\bsnm{Maynard~Smith},~\bfnm{J.}\binits{J.}}
(\byear{1978}).
\btitle{The Evolution of Sex}.
\bpublisher{Cambridge Univ. Press}, \baddress{Cambridge}.
\bptok{imsref}%
\end{bbook}
\endbibitem

%b29 ###
\bibitem[\protect\citeauthoryear{Muller}{1964}]{Muller1964}
\begin{barticle}[pbm]
\bauthor{\bsnm{Muller},~\bfnm{H.~J.}\binits{H.~J.}}
(\byear{1964}).
\btitle{The relation of recombination to mutational advance}.
\bjournal{Mutat. Res.}
\bvolume{106}
\bpages{2--9}.
\bid{issn={0027-5107}, pmid={14195748}}
\bptok{imsref}%
\end{barticle}
\endbibitem

%b30 ###
\bibitem[\protect\citeauthoryear{Park and Krug}{2007}]{ParkKrug2007}
\begin{barticle}[auto:STB|2012/04/25|11:03:30]
\bauthor{\bsnm{Park},~\bfnm{S.~C.}\binits{S.~C.}} \AND
  \bauthor{\bsnm{Krug},~\bfnm{J.}\binits{J.}}
(\byear{2007}).
\btitle{Clonal interference in large populations}.
\bjournal{PNAS}
\bvolume{104}
\bpages{18135--18140}.
\bptok{imsref}%
\end{barticle}
\endbibitem

%b31 ###
\bibitem[\protect\citeauthoryear{Poon and Chao}{2005}]{PoonChao2005}
\begin{barticle}[pbm]
\bauthor{\bsnm{Poon},~\bfnm{Art}\binits{A.}} \AND
  \bauthor{\bsnm{Chao},~\bfnm{Lin}\binits{L.}}
(\byear{2005}).
\btitle{The rate of compensatory mutation in the DNA bacteriophage phiX174}.
\bjournal{Genetics}
\bvolume{170}
\bpages{989--999}.
\bid{doi={10.1534/genetics.104.039438}, issn={0016-6731},
  pii={genetics.104.039438}, pmcid={1451187}, pmid={15911582}}
\bptok{imsref}%
\end{barticle}
\endbibitem

%b32 ###
\bibitem[\protect\citeauthoryear{Protter}{2004}]{Protter2004}
\begin{bbook}[mr]
\bauthor{\bsnm{Protter},~\bfnm{Philip~E.}\binits{P.~E.}}
(\byear{2004}).
\btitle{Stochastic Integration and Differential Equations},
\bedition{2nd} ed.
\bseries{Applications of Mathematics (New York)}
\bvolume{21}.
\bpublisher{Springer}, \baddress{Berlin}.
\bid{mr={2020294}}
\bptok{imsref}%
\end{bbook}
\endbibitem

%b33 ###
\bibitem[\protect\citeauthoryear{Rouzine, Wakeley and
  Coffin}{2003}]{RouzineWakeleyCoffin2003}
\begin{barticle}[pbm]
\bauthor{\bsnm{Rouzine},~\bfnm{Igor~M.}\binits{I.~M.}},
  \bauthor{\bsnm{Wakeley},~\bfnm{John}\binits{J.}} \AND
  \bauthor{\bsnm{Coffin},~\bfnm{John~M.}\binits{J.~M.}}
(\byear{2003}).
\btitle{The solitary wave of asexual evolution}.
\bjournal{Proc. Natl. Acad. Sci. USA}
\bvolume{100}
\bpages{587--592}.
\bid{doi={10.1073/pnas.242719299}, issn={0027-8424}, pii={242719299},
  pmcid={141040}, pmid={12525686}}
\bptok{imsref}%
\end{barticle}
\endbibitem

%b34 ###
\bibitem[\protect\citeauthoryear{Shiga and Shimizu}{1980}]{ShigaShimizu1980}
\begin{barticle}[mr]
\bauthor{\bsnm{Shiga},~\bfnm{Tokuzo}\binits{T.}} \AND
  \bauthor{\bsnm{Shimizu},~\bfnm{Akinobu}\binits{A.}}
(\byear{1980}).
\btitle{Infinite-dimensional stochastic differential equations and their
  applications}.
\bjournal{J. Math. Kyoto Univ.}
\bvolume{20}
\bpages{395--416}.
\bid{issn={0023-608X}, mr={0591802}}
\bptok{imsref}%
\end{barticle}
\endbibitem

%b35 ###
\bibitem[\protect\citeauthoryear{Stephan, Chao and
  Smale}{1993}]{StephanChaoSmale1993}
\begin{barticle}[auto:STB|2012/04/25|11:03:30]
\bauthor{\bsnm{Stephan},~\bfnm{W.}\binits{W.}},
  \bauthor{\bsnm{Chao},~\bfnm{L.}\binits{L.}} \AND
  \bauthor{\bsnm{Smale},~\bfnm{J.}\binits{J.}}
(\byear{1993}).
\btitle{The advance of Muller's ratchet in a~haploid asexual population:
  Approximate solutions based on diffusion theory}.
\bjournal{Genet. Res.}
\bvolume{61}
\bpages{225--231}.
\bptok{imsref}%
\end{barticle}
\endbibitem

%b36 ###
\bibitem[\protect\citeauthoryear{Wagner and Gabriel}{1990}]{WagnerGabriel1990}
\begin{barticle}[auto:STB|2012/04/25|11:03:30]
\bauthor{\bsnm{Wagner},~\bfnm{G.}\binits{G.}} \AND
  \bauthor{\bsnm{Gabriel},~\bfnm{W.}\binits{W.}}
(\byear{1990}).
\btitle{What stops Muller's ratchet in the absence of recombination?}
\bjournal{Evolution}
\bvolume{44}
\bpages{715--731}.
\bptok{imsref}%
\end{barticle}
\endbibitem

%b37 ###
\bibitem[\protect\citeauthoryear{Waxman and Loewe}{2010}]{WaxmanLoewe2010}
\begin{barticle}[auto:STB|2012/04/25|11:03:30]
\bauthor{\bsnm{Waxman},~\bfnm{D.}\binits{D.}} \AND
  \bauthor{\bsnm{Loewe},~\bfnm{L.}\binits{L.}}
(\byear{2010}).
\btitle{A stochastic model for a single click of Muller's ratchet}.
\bjournal{J.~Theor. Biol.}
\bvolume{264}
\bpages{1120--1132}.
\bptok{imsref}%
\end{barticle}
\endbibitem

%b38 ###
\bibitem[\protect\citeauthoryear{Yu, Etheridge and
  Cuthbertson}{2010}]{YuEtheridgeCuthbertson2010}
\begin{barticle}[mr]
\bauthor{\bsnm{Yu},~\bfnm{Feng}\binits{F.}},
  \bauthor{\bsnm{Etheridge},~\bfnm{Alison}\binits{A.}} \AND
  \bauthor{\bsnm{Cuthbertson},~\bfnm{Charles}\binits{C.}}
(\byear{2010}).
\btitle{Asymptotic behavior of the rate of adaptation}.
\bjournal{Ann. Appl. Probab.}
\bvolume{20}
\bpages{978--1004}.
\bid{doi={10.1214/09-AAP645}, issn={1050-5164}, mr={2680555}}
\bptok{imsref}%
\end{barticle}
\endbibitem

\end{thebibliography}
\end{document}